\documentclass[11pt,a4paper]{article}

\usepackage[utf8]{inputenc}
\PassOptionsToPackage{hyphens}{url}
\usepackage{amsmath,amssymb,amsthm}
\usepackage{mathtools}
\usepackage{graphicx}
\usepackage{booktabs}
\usepackage{array}
\usepackage[hypertexnames=false,colorlinks=true,linkcolor=blue,citecolor=blue,urlcolor=blue,pdfauthor={Kyle Elliott Mathewson},pdftitle={Balance in Inequality: Testing the Physical Equilibrium, Precession, and Structural Stability of Magic Gems}]{hyperref}
\usepackage{geometry}
\usepackage{float}

\makeatletter
\def\Cref@getprefix#1:#2\relax{#1}
\def\Cref@prefix@sec{sec}
\def\Cref@prefix@subsec{subsec}
\def\Cref@prefix@fig{fig}
\def\Cref@prefix@tab{tab}
\def\Cref@prefix@eq{eq}
\def\Cref@prefix@thm{thm}
\def\Cref@prefix@prop{prop}
\def\Cref@prefix@lem{lem}
\def\Cref@prefix@cor{cor}
\def\Cref@prefix@conj{conj}
\def\Cref@prefix@def{def}
\def\Cref@prefix@rem{rem}
\def\Cref@prefix@ex{ex}

\def\Cref@type@cap#1{%
  \edef\Cref@tmp{#1}%
  \ifx\Cref@tmp\Cref@prefix@sec Section%
  \else\ifx\Cref@tmp\Cref@prefix@subsec Section%
  \else\ifx\Cref@tmp\Cref@prefix@fig Figure%
  \else\ifx\Cref@tmp\Cref@prefix@tab Table%
  \else\ifx\Cref@tmp\Cref@prefix@eq Equation%
  \else\ifx\Cref@tmp\Cref@prefix@thm Theorem%
  \else\ifx\Cref@tmp\Cref@prefix@prop Proposition%
  \else\ifx\Cref@tmp\Cref@prefix@lem Lemma%
  \else\ifx\Cref@tmp\Cref@prefix@cor Corollary%
  \else\ifx\Cref@tmp\Cref@prefix@conj Conjecture%
  \else\ifx\Cref@tmp\Cref@prefix@def Definition%
  \else\ifx\Cref@tmp\Cref@prefix@rem Remark%
  \else\ifx\Cref@tmp\Cref@prefix@ex Example%
  \else Ref.%
  \fi\fi\fi\fi\fi\fi\fi\fi\fi\fi\fi\fi\fi
}

\def\Cref@type@low#1{%
  \edef\Cref@tmp{#1}%
  \ifx\Cref@tmp\Cref@prefix@sec section%
  \else\ifx\Cref@tmp\Cref@prefix@subsec section%
  \else\ifx\Cref@tmp\Cref@prefix@fig figure%
  \else\ifx\Cref@tmp\Cref@prefix@tab table%
  \else\ifx\Cref@tmp\Cref@prefix@eq equation%
  \else\ifx\Cref@tmp\Cref@prefix@thm theorem%
  \else\ifx\Cref@tmp\Cref@prefix@prop proposition%
  \else\ifx\Cref@tmp\Cref@prefix@lem lemma%
  \else\ifx\Cref@tmp\Cref@prefix@cor corollary%
  \else\ifx\Cref@tmp\Cref@prefix@conj conjecture%
  \else\ifx\Cref@tmp\Cref@prefix@def definition%
  \else\ifx\Cref@tmp\Cref@prefix@rem remark%
  \else\ifx\Cref@tmp\Cref@prefix@ex example%
  \else ref.%
  \fi\fi\fi\fi\fi\fi\fi\fi\fi\fi\fi\fi\fi
}

\DeclareRobustCommand{\Cref}[1]{%
  \begingroup
  \edef\Cref@label{#1}%
  \edef\Cref@prefix{\expandafter\Cref@getprefix\Cref@label:\relax}%
  \hyperref[#1]{\Cref@type@cap{\Cref@prefix}~\ref*{#1}}%
  \endgroup
}
\DeclareRobustCommand{\cref}[1]{%
  \begingroup
  \edef\Cref@label{#1}%
  \edef\Cref@prefix{\expandafter\Cref@getprefix\Cref@label:\relax}%
  \hyperref[#1]{\Cref@type@low{\Cref@prefix}~\ref*{#1}}%
  \endgroup
}
\makeatother

\makeatletter
\long\def\@makecaption#1#2{\vskip\abovecaptionskip\small\sbox\@tempboxa{#1: #2}\ifdim\wd\@tempboxa>\hsize #1: #2\par\else\global\@minipagefalse\hb@xt@\hsize{\hfil\box\@tempboxa\hfil}\fi\vskip\belowcaptionskip}
\makeatother

\theoremstyle{plain}
\newtheorem{theorem}{Theorem}[section]

\newtheorem{proposition}[theorem]{Proposition}
\newtheorem{corollary}[theorem]{Corollary}

\theoremstyle{definition}
\newtheorem{definition}[theorem]{Definition}

\theoremstyle{remark}
\newtheorem{remark}[theorem]{Remark}

\newcommand{\R}{\mathbb{R}}
\newcommand{\Z}{\mathbb{Z}}

\newcommand{\Q}{\mathbb{Q}}
\newcommand{\tr}{\operatorname{tr}}
\newcommand{\conv}{\operatorname{conv}}
\newcommand{\vertset}{\operatorname{vert}}
\newcommand{\vol}{\operatorname{vol}}
\newcommand{\cone}{\operatorname{cone}}
\newcommand{\diag}{\operatorname{diag}}

\title{\textbf{Balance in Inequality: Testing the Physical Equilibrium,\\ Precession, and Structural Stability of Magic Gems}}

\author{Kyle Elliott Mathewson\\
\textit{Faculty of Science, University of Alberta}\\
\texttt{kyle.mathewson@ualberta.ca}}

\date{}

\begin{document}

\maketitle

\begin{abstract}
A magic square is a combinatorial paradigm of \emph{balance in inequality}: its $n^2$ entries $\{1,\dots,n^2\}$ are strictly distinct, yet every row, column, and main diagonal sums to the same constant $\mu_n=n(n^2+1)/2$. Lifting the centred grid cells into $\R^3$ as radial vectors $\mathbf p_{ij}(\lambda)=(x_j,y_i,\lambda z_{ij})$ and taking their convex hull—the \emph{Magic Gem} $\mathcal G_\lambda(S)$—raises a central mechanical question: when does combinatorial line-sum equality suffice for physical equilibrium in three dimensions, and when a nonlinear physical realization breaks that balance, how can an opposing inequality restore it? We test this principle across all $880$ order-four magic squares by progressing from one-dimensional radial spokes to three-dimensional rotation, continuous solid and skeletal hulls, and freestanding stacked columns. For radial spokes from the origin $O$, natural vector lengths balance identically, whereas trimming the spokes to equal unit length breaks equilibrium on all $832$ non-associative squares in 3D and on every square in 2D, where the irrational ring-radius ratio $r_2/r_1=\sqrt{5}$ bounds the resultant below by the Pell convergent $9/4\approx\sqrt{5}$ ($18\sqrt{10}/5-8\sqrt{2}\approx0.0705$). By length--weight duality, radial equilibrium is restored either by assigning positive unequal tip weights $\mathbf w^*>\mathbf 0$ on the unit sphere $S^2$ or by scaling spoke lengths as a power law ($r^{\alpha^*}$ in 2D and cubic vertical vanishing $(\mathbf V_3)_z=0$ in 3D). In three-dimensional rotation, equal point masses at the unequal vertices $\mathbf p_{ij}(\lambda)$ make every semi-magic square an axisymmetric top with universal inertia ratio $I_\perp(\lambda)/I_{zz}=\frac12(1+\lambda^2(n^2+1))$ and spherical-top transition $\lambda_{\mathrm{iso}}=(n^2+1)^{-1/2}$, whereas weighting the vertices by their unequal entries $m_{ij}=s_{ij}$ exposes the quadrant quadrupole $Q_{xy}=\sum z_{ij}x_jy_i$, splitting the horizontal principal moments by $2|Q_{xy}|$ and triggering $180^\circ$ Dzhanibekov flips on $824$ squares. Passing from discrete points to the uniform solid hull $\mathcal G_0(S)$ shifts the centre of mass off $O$ on $829$ squares because interior grid points $P_S^{(1)}$ are omitted and hull facets are weighted by unequal simplex volumes; because the outer hull vertices and inner core points pull in strictly opposite directions ($V_0\overline{\mathbf p}_{\mathrm{outer}}+V_1\overline{\mathbf p}_{\mathrm{inner}}=\mathbf 0$), re-engaging the inner core $\mathcal G_1(S)=\conv P_S^{(1)}$ restores balance via two-density core inclusions, two-gauge nested wireframes, or the complete $\binom{n^2}{2}$-rod network $K_{n^2}$ ($\Sigma_{\mathrm{rods}}=\frac{2n^2-3}{3(n^2-1)}\Sigma(\lambda)$), while every boundary hull unfolds into a non-overlapping papercraft net with $f_0-1$ glue tabs. Finally, compressing the vertical aspect to the cubic print scale $\lambda=1/5$ stabilizes $74.5\%$ of all gem facets on a flat table, and alternating any gem $\mathcal G_\lambda(S)$ with its inverted mirror twin $\mathcal G_\lambda(r^2\overline S)=-\mathcal G_\lambda(S)$ cancels opposing facet tilts to lock chiral blocks into a straight vertical tower.

\smallskip\noindent\textbf{Keywords:} magic squares; convex polyhedra; Albrecht D\"urer; centre of mass; rigid-body precession; Dzhanibekov effect; polyhedral nets; mathematical sculpture.\\
\noindent\textbf{MSC 2020:} 52B10, 05B15, 70E15, 00A66.
\end{abstract}

\section{Introduction: Balance in Inequality}
\label{sec:intro}

Over five centuries ago, Albrecht D\"urer contributed two enduring landmarks at the intersection of mathematics and the visual arts: the order-four magic square in his 1514 engraving \emph{Melencolia~I} (square \#175 in Fr\'enicle's 1693 census of the $880$ order-four magic squares \cite{frenicle1693,dudeney1917amusements,andrews1960magic}) and the first printed polyhedral net unfoldings in his 1525 treatise \emph{Underweysung der Messung} \cite{durer1525underweysung,shephard1975convex,demaine2007geometric,ghomi2014durer}. A magic square is a quintessential study in \emph{balance in inequality}: an arrangement $S=(s_{ij})_{0\le i,j<n}$ of the first $N=n^2$ positive integers $\{1,\dots,N\}$ in which no two cells are equal, yet every row, column, and main diagonal sums to the exact same magic constant $\mu_n=n(N+1)/2$ (a \emph{semi-magic square} if the row and column sums equal $\mu_n$). Measuring coordinates from the board centre and the mean entry $\bar s=(N+1)/2$,
\begin{equation}
\label{eq:coords}
x_j = j - \frac{n-1}{2},\qquad
y_i = \frac{n-1}{2} - i,\qquad
z_{ij} = s_{ij} - \frac{N+1}{2},
\end{equation}
lifts each cell $(i,j)$ into a three-dimensional point and radial vector from the origin $O=(0,0,0)$:
\[
\mathbf p_{ij}(\lambda) = (x_j,\, y_i,\, \lambda z_{ij}) \in \R^3,\qquad
P_S(\lambda) = \bigl\{\mathbf p_{ij}(\lambda) : 0\le i,j<n\bigr\},
\]
where $\lambda>0$ is the vertical lift scale (\Cref{fig:origin-vectors-3d}(a)). The convex hull
\[
\mathcal G_\lambda(S) := \conv P_S(\lambda) \subset \R^3
\qquad\bigl(\mathcal G(S):=\mathcal G_1(S)\bigr)
\]
is a three-dimensional convex polyhedron that we call the \emph{Magic Gem} of $S$, reuniting D\"urer's magic squares with his polyhedral geometry.

In the flat grid plane, Loly \cite{loly2004invariance} and Rogers and Loly \cite{rogers2004inertia,rogers2005multipole} showed that interpreting the unequal entries $s_{ij}$ as point masses at $\mathbf r_{ij}=(x_j,y_i)$ places the planar centre of mass at the origin $O$ ($\sum_{i,j}z_{ij}\mathbf r_{ij}=\mathbf 0$) and makes the polar moment of inertia $\sum_{i,j}s_{ij}(x_j^2+y_i^2)=n^2(n^4-1)/12$ invariant across all semi-magic squares of order $n$. When a lifted magic square is fabricated as a three-dimensional physical object—whether as a radial spoke star, a planar balance wheel, a spinning top on a pivot, a cast solid, a 3D-printed wireframe, a folded papercraft shell, or a freestanding tower of blocks \cite{mathewson2025magicgemweb}—its unequal entries are subjected to progressively nonlinear geometric transformations. This motivates our central hypothesis: \emph{combinatorial line-sum equality guarantees physical equilibrium whenever a mechanical observable is linear or symmetric-quadratic in the grid coordinates, neutralizing entry inequality universally across all $880$ order-four magic squares; when a nonlinear geometric constraint breaks that equilibrium, balance can be restored by pairing two opposing inequalities that cancel each other.}

\begin{figure}[htbp]
\centering
\includegraphics[width=\textwidth]{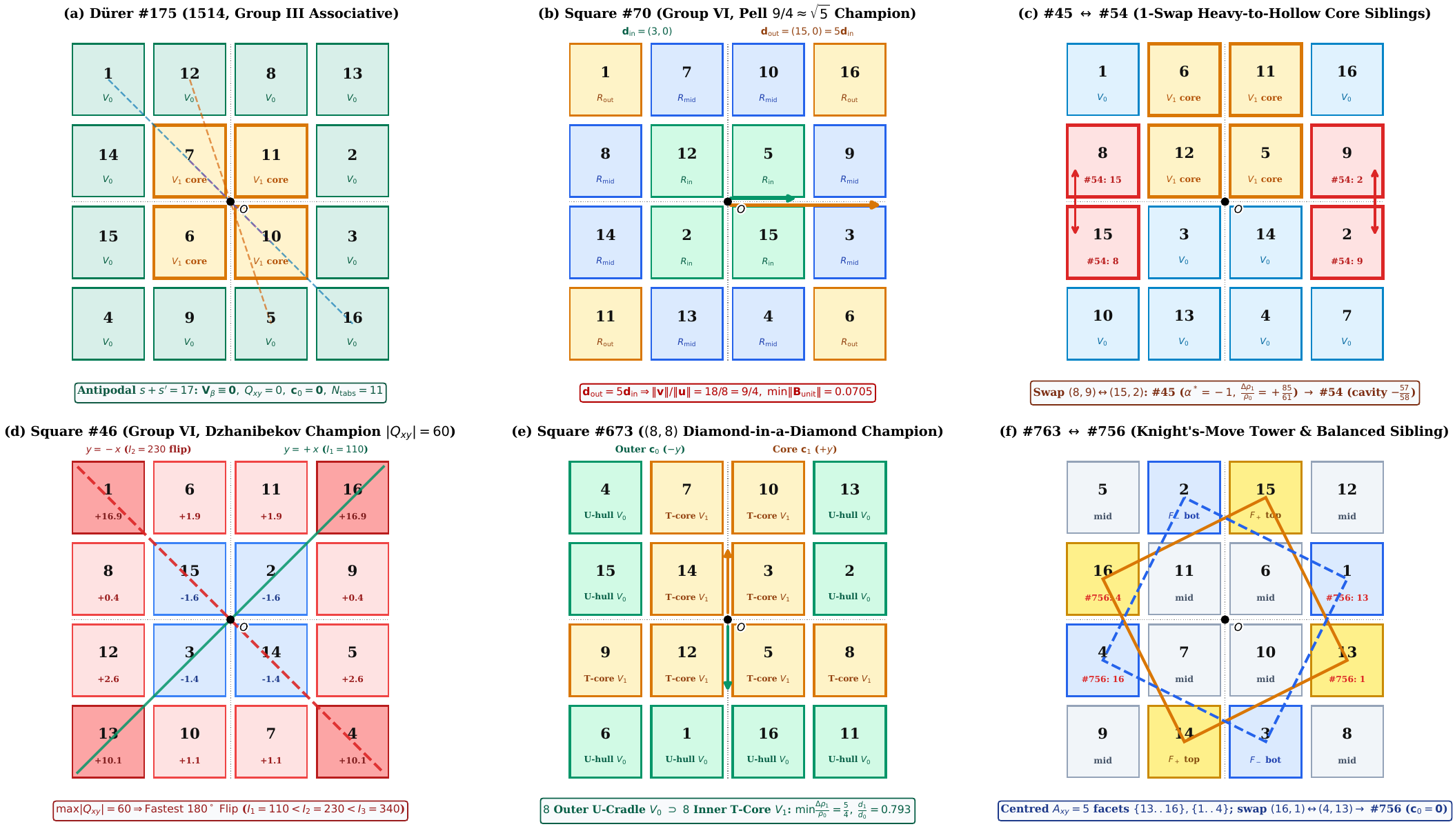}
\caption{Visual anatomy of the six flagship order-four magic squares examined across Sections~\ref{sec:vectors}--\ref{sec:stacking}. (a)~D\"urer's associative square \#175 (1514, Group~III), whose antipodal pairs sum to $17$ across $O$, partitioning into $12$ outer hull vertices $V_0$ and $4$ inner core vertices $V_1$. (b)~Square \#70 (Group~VI), coloured by concentric rings $R_{\mathrm{in}},R_{\mathrm{mid}},R_{\mathrm{out}}$, where parallel ring dipoles $\mathbf d_{\mathrm{out}}=5\mathbf d_{\mathrm{in}}$ achieve the Pell ratio $9/4\approx\sqrt{5}$, and (c)~square \#45 ($\kappa=+9$, $\alpha^*=-1$), whose single border swap $(8,9)\leftrightarrow(15,2)$ flips a heavy-core gem into hollow-cavity sibling \#54. (d)~Square \#46 (Group~VI), shaded by quadrant quadrupole $z_{ij}x_jy_i$, attaining the maximum $|Q_{xy}|=60$ and fastest $180^\circ$ Dzhanibekov flip about $y=-x$. (e)~Square \#673 (Group~VI), whose $8$-vertex U-cradle $V_0$ encloses an $8$-vertex inverted-T core $V_1$ ($\Delta\rho_1/\rho_0=5/4$), and (f)~square \#763 (Group~VI), whose centred knight's-move facets $\{13,14,15,16\}$ and $\{1,2,3,4\}$ lock into an area-$5$ mirror-twin tower.}
\label{fig:hero-squares}
\end{figure}

We test this hypothesis across four progressive mechanical regimes, illustrated on six recurring flagship squares (\Cref{fig:hero-squares}). In \Cref{sec:vectors}, we examine one-dimensional radial spoke sculptures radiating from the origin $O$. At natural spoke lengths, every magic square balances at $O$ because unequal lengths compensate unequal directions. Trimming all sixteen spokes to equal unit length breaks equilibrium on all $832$ non-associative squares in 3D and on \emph{every} square in 2D, where the middle grid ring sits at an irrational radius ratio $r_2/r_1=\sqrt{5}$ relative to the inner and outer rings ($r_3/r_1=3$) and prevents exact cancellation beyond the Pell convergent $9/4\approx\sqrt{5}$ (square \#70, \Cref{fig:hero-squares}(b), \Cref{thm:three-ring}). By length--weight duality (\Cref{prop:reciprocity}), equilibrium is restored either by playing aligned inner and outer ring dipoles against the middle ring via a power-law length profile $r^{\alpha^*}$ (such as $\alpha^*=-1$ on \#45, \Cref{fig:hero-squares}(c), \Cref{cor:planar-power-roots}), or by assigning positive unequal tip weights $\mathbf w^*>\mathbf 0$ on the unit sphere $S^2$ and cubic spoke lengths ($\beta=3$, where the vertical resultant vanishes universally, \Cref{thm:3d-power}).

In \Cref{sec:dynamics}, we set balanced sculptures into three-dimensional rotation on a central pivot under gravity. When the sixteen unequal vertices $\mathbf p_{ij}(\lambda)$ carry equal masses, line-sum equality eliminates all off-diagonal products of inertia, turning every semi-magic square into a smooth axisymmetric top that passes through an exact spherical-top transition at $\lambda_{\mathrm{iso}}=(n^2+1)^{-1/2}$ ($1/\sqrt{17}$ at $n=4$, \Cref{thm:rotor}). Replacing equal masses with the unequal magic entries $m_{ij}=s_{ij}$ (\Cref{thm:quadrupole-split}) breaks this rotational symmetry by exposing the planar quadrupole moment $Q_{xy}=\sum z_{ij}x_jy_i$ and the cubic vector resultant from \Cref{sec:vectors}, splitting the horizontal principal axes by $2|Q_{xy}|$ and triggering $180^\circ$ Dzhanibekov flips on $824$ of the $880$ squares (led by \#46 at $|Q_{xy}|=60$, \Cref{fig:hero-squares}(d)).

In \Cref{sec:nested}, we pass from discrete point clouds to continuous polyhedral hulls $\mathcal G_0(S)=\conv P_S$. Casting $\mathcal G_0(S)$ as a uniform solid preserves near-axisymmetric inertia about its own centroid $\mathbf c_0(S)$, but shifts $\mathbf c_0(S)$ off the pivot $O$ on $829$ of the $880$ squares ($51$ balance naturally, \Cref{prop:51-balanced}) because interior grid points $P_S^{(1)}$ are omitted and boundary facets are weighted by unequal simplex volumes. Because $\sum\mathbf p_{ij}=\mathbf 0$, the outer hull vertices and the inner core points pull in strictly opposite directions across the origin ($V_0\overline{\mathbf p}_{\mathrm{outer}}=-V_1\overline{\mathbf p}_{\mathrm{inner}}$, \Cref{prop:opposition}). Physically re-engaging the inner core $\mathcal G_1(S)=\conv P_S^{(1)}$ restores $\mathbf c=\mathbf 0$ across composite solids (such as denser-core gem \#673 and hollow-cavity gem \#54, \Cref{fig:hero-squares}(c,e), \Cref{thm:diamond-balance}), single-material two-gauge wireframes and complete $K_{16}$ strut networks (\Cref{thm:complete-network}), and unfolded papercraft shells with $f_0-1$ glue tabs (\Cref{prop:papercraft-net}).

Finally, in \Cref{sec:stacking}, we remove the central pivot and place Magic Gems onto a flat horizontal table under gravity. While tall $\lambda=1$ gems topple on more than half of their facets, compressing the vertical lift to the cubic print scale $\lambda_{\mathrm{cube}}=1/5$ stabilizes $74.5\%$ of all $14{,}824$ facets (\Cref{prop:stability-census}). When stacking multiple blocks vertically, identical copies of a single gem wedge apart whenever the top and bottom facet normals fail to oppose ($\mathbf n_-\ne-\mathbf n_+$ on $756$ gems); alternating $\mathcal G_\lambda(S)$ with its inverted mirror twin $\mathcal G_\lambda(r^2\overline S)=-\mathcal G_\lambda(S)$ (\Cref{thm:mirror-tower}) pairs every slanted facet with an equal and opposite partner, locking chiral knight's-move blocks (\#763, \Cref{fig:hero-squares}(f)) into a straight vertical column (\Cref{prop:knights-move}).

Across all four regimes, the order-four censuses stratify into three symmetry tiers (\Cref{sec:discussion}): universal properties governed by linear line sums hold across all $880$ squares; unweighted nonlinear cancellations hold on the $48$ associative squares of Dudeney's Group~III ($r^2S=\overline S$, led by D\"urer's \#175), where central inversion $\mathbf p_{-c}=-\mathbf p_c$ cancels opposite cells pairwise; and one-parameter counterbalancing between opposing rings or hull layers occurs in the $304$ self-complementary squares of Dudeney's Group~VI ($\overline S=hS$ or $vS$), where a $180^\circ$ half-turn symmetry ($C_2\subset SO(3)$) confines residual imbalance vectors to a single coordinate axis. Every theorem and census table is certified in exact rational arithmetic by \texttt{anc/verify\_paper\_c.py} (\Cref{sec:appendix}) and accompanied by real-time 3D simulations in \texttt{web/index.html}.

Throughout the paper, $A_n=\{(x_j,y_i):0\le i,j<n\}$ denotes the centred grid \eqref{eq:coords}, $a_n=(n^2-1)/12$ and $b_n=(n^4-1)/12=(n^2+1)a_n$ denote the marginal variances of $x_j$ and $z_{ij}$, and the dihedral group $D_4=\langle r,h\rangle$ (with $90^\circ$ rotation $r$ and reflections $h,v$) acts on arrangements together with the complement involution $\tau(S)=\overline S=(N+1)\mathbf 1-S$. Companion preprints \cite{mathewson2026gems,mathewson2026moments} treat the static polyhedral combinatorics and bivariate Prouhet--Tarry--Escott moment hierarchy of semi-magic squares; the present paper is self-contained. Order-four magic squares are numbered $\#1,\dots,\#880$ by their 1-based position in Fr\'enicle's standard list (\texttt{anc/squares880.txt}).

\section{Radial Spokes from the Origin and Length--Weight Duality}
\label{sec:vectors}

In the fabrication suite \texttt{magicGemWeb} \cite{mathewson2025magicgemweb}, the sixteen cells of an order-four magic square are built as two radial sculptures pivoted at the origin $O$: a three-dimensional sixteen-spoke hub-and-strut star (\texttt{magicGem.scad}, \Cref{fig:origin-vectors-3d}(a)) with spokes $\mathbf p_{ij}(\lambda)=(x_j,y_i,\lambda z_{ij})$ of length $R_{ij}(\lambda)=\|\mathbf p_{ij}(\lambda)\|_2$ along unit directions $\hat{\mathbf p}_{ij}(\lambda)=\mathbf p_{ij}(\lambda)/R_{ij}(\lambda)\in S^2$, and a planar sixteen-spoke balance wheel on a conical fulcrum (\texttt{weighted\_ring.stl}, \Cref{fig:ring-vectors}(a)) with spokes along $\mathbf r_{ij}=(x_j,y_i)$ of radius $r_{ij}=\|\mathbf r_{ij}\|_2$ and unit direction $\mathbf u_{ij}=\mathbf r_{ij}/r_{ij}\in S^1$. At natural linear lengths ($L_{ij}=R_{ij}(\lambda)$ or $\ell(r_{ij})=r_{ij}$), both sculptures balance identically at $O$ because $\sum_{i,j}\mathbf p_{ij}(\lambda)=\mathbf 0$ and $\sum_{i,j}z_{ij}\mathbf r_{ij}=\mathbf 0$ \cite{loly2004invariance,rogers2005multipole}, closing the sixteen planar vectors into a head-to-tail polygon (\Cref{fig:ring-vectors}(b)).

Our first test probes the mechanism of this equilibrium: does radial balance reside in the directions of the spokes alone, or does it rely on the unequal spoke lengths compensating the unequal directions? And when equalizing spoke lengths breaks balance, how can an opposing inequality in tip weights or power-law lengths restore it?

\subsection{Length--Weight Duality: Why Equalizing Spoke Lengths Breaks Balance}

\begin{definition}[Generalized 3D and planar origin-vector resultants]
\label{def:radial-resultant}
For an arrangement $S$ of order $n$ and vertical lift scale $\lambda>0$, let $R_{ij}(\lambda)=\sqrt{x_j^2+y_i^2+\lambda^2z_{ij}^2}$ and $\hat{\mathbf p}_{ij}(\lambda)=\mathbf p_{ij}(\lambda)/R_{ij}(\lambda)$ ($\mathbf 0$ if $R_{ij}(\lambda)=0$). Given cell weights $w=(w_{ij})$ and prescribed 3D spoke lengths $L=(L_{ij})$, the \emph{3D origin-vector resultant} is
\begin{equation}
\label{eq:V-3d-general}
\mathbf V_{w,L}(S;\lambda) = \sum_{i,j} w_{ij}\,L_{ij}\,\hat{\mathbf p}_{ij}(\lambda)
= \sum_{R_{ij}(\lambda)>0} \left(\frac{w_{ij}L_{ij}}{R_{ij}(\lambda)}\right)\mathbf p_{ij}(\lambda) \in \R^3.
\end{equation}
For equal weights $w_{ij}=1$ and power-law 3D lengths $L_{ij}=R_{ij}(\lambda)^\beta$ ($\beta\in\R$), we write
\begin{equation}
\label{eq:V-beta}
\mathbf V_\beta(S;\lambda) := \sum_{i,j} R_{ij}(\lambda)^\beta\,\hat{\mathbf p}_{ij}(\lambda)
= \sum_{i,j} R_{ij}(\lambda)^{\beta-1}\,\mathbf p_{ij}(\lambda) \in \R^3,
\end{equation}
and for weighted unit vectors ($L_{ij}=1$) on $S^2$, we write $\mathbf B^{(3)}_{w,\lambda}(S):=\mathbf V_{w,1}(S;\lambda)=\sum_{i,j}w_{ij}\hat{\mathbf p}_{ij}(\lambda)$. Similarly, in the grid plane, given cell masses $m=(m_{ij})$ and a radial length profile $\ell\colon[0,\infty)\to[0,\infty)$ along unit rays $\mathbf u_{ij}=\mathbf r_{ij}/r_{ij}$, the \emph{planar radial resultant} is
\begin{equation}
\label{eq:B-planar}
\mathbf B_{m,\ell}(S) = \sum_{i,j} m_{ij}\,\ell(r_{ij})\,\mathbf u_{ij}
= \sum_{r_{ij}>0} \left(\frac{m_{ij}\ell(r_{ij})}{r_{ij}}\right)(x_j,y_i) \in \R^2.
\end{equation}
\end{definition}

Because $\ell(r_{ij})$ depends only on the radius $r_{ij}$, grid symmetry gives $\sum_{i,j}\ell(r_{ij})\mathbf u_{ij}=\mathbf 0$, so replacing the positive entry masses $m_{ij}=s_{ij}=\bar s+z_{ij}$ by the signed heights $m_{ij}=z_{ij}$ leaves the planar resultant unchanged: $\mathbf B_{s,\ell}(S)=\mathbf B_{z,\ell}(S)$.

\begin{proposition}[Length--weight duality and the order-three two-ring law]
\label{prop:reciprocity}
Let $n\ge3$ and $\lambda>0$. Scaling 3D vector lengths from the origin to $L_{ij}$ with weights $w_{ij}$ produces the same resultant \eqref{eq:V-3d-general} as keeping the natural 3D vectors $\mathbf p_{ij}(\lambda)$ ($L'_{ij}=R_{ij}(\lambda)$) and assigning effective weights $w'_{ij}=w_{ij}L_{ij}/R_{ij}(\lambda)$. In particular, power-law 3D spoke lengths $L_{ij}=R_{ij}(\lambda)^\beta$ with unit weights ($w_{ij}=1$) are mechanically dual to natural 3D vectors $\mathbf p_{ij}(\lambda)$ with weights $w'_{ij}=R_{ij}(\lambda)^{\beta-1}$, and unit-length 3D vectors on $S^2$ ($L_{ij}=1$) balance identically ($\sum_{i,j}w_{ij}\hat{\mathbf p}_{ij}(\lambda)=\mathbf 0$) for every arrangement $S$ when weighted by their unequal 3D lengths $w_{ij}=c\,R_{ij}(\lambda)$ ($c\in\R$). In the grid plane, unit vectors $\mathbf u_{ij}$ balance for every semi-magic square ($\sum_{i,j}m_{ij}\mathbf u_{ij}=\mathbf 0$) whenever $m_{ij}=c\,z_{ij}r_{ij}$, and a radial profile $\ell_{ij}$ with $m_{ij}=z_{ij}$ balances every semi-magic square if and only if $(\ell_{ij}/r_{ij})\mathbf r_{ij}$ has components in the row--column subspace $W_{\mathrm{rc}}=\{a(x_j)+b(y_i)\}$.

At order $n=3$, where the eight non-central cells of $A_3$ lie on two concentric rings—four edge midpoints $R_1=\{(\pm1,0),(0,\pm1)\}$ ($r_1=1$) and four corners $R_{\sqrt{2}}=\{(\pm1,\pm1)\}$ ($r_2=\sqrt{2}$)—every semi-magic square $S$ satisfies $\mathbf d_1(S)+\mathbf d_2(S)=\mathbf 0$ for $\mathbf d_k(S)=\sum_{(i,j)\in R_k}z_{ij}\mathbf r_{ij}$, and hence
\begin{equation}
\label{eq:n3-radial}
\mathbf B_{z,\ell}(S) = \left(\ell(1) - \frac{\ell(\sqrt{2})}{\sqrt{2}}\right)\mathbf d_1(S).
\end{equation}
For all eight magic squares of order~$3$ (the $D_4$-orbit of the Lo Shu square), $\|\mathbf d_1(S)\|_2=\sqrt{80}=4\sqrt{5}$, giving
\begin{equation}
\label{eq:n3-norm}
\|\mathbf B_{z,\ell}(S)\|_2 = 4\sqrt{5}\,\left|\ell(1) - \frac{\ell(\sqrt{2})}{\sqrt{2}}\right|,
\qquad
\|\mathbf B_{\mathrm{unit}}(S)\|_2 = 4\sqrt{5} - 2\sqrt{10} \approx 2.619717,
\end{equation}
so $\mathbf B_{z,\ell}(S)=\mathbf 0$ at $n=3$ if and only if $\ell(\sqrt{2})=\sqrt{2}\,\ell(1)$.
\end{proposition}

\begin{proof}
The duality relations follow directly from \eqref{eq:V-3d-general}--\eqref{eq:B-planar}, the identity $\sum_{i,j}\mathbf p_{ij}(\lambda)=\mathbf 0$, and the characterization of $W_{\mathrm{rc}}$ as the orthogonal complement of the zero-line-sum space \cite[Sect.~2]{mathewson2026moments}. At order $n=3$, $\mathbf r_{11}=(0,0)$ implies $\mathbf d_1(S)+\mathbf d_2(S)=\sum_{i,j}z_{ij}\mathbf r_{ij}=\mathbf 0$ for every semi-magic square, which yields \eqref{eq:n3-radial} upon grouping \eqref{eq:B-planar} over $R_1$ and $R_{\sqrt{2}}$. Evaluating $\mathbf d_1=(-8,4)$ on the Lo Shu square gives $\|\mathbf d_1\|_2=\sqrt{80}=4\sqrt{5}$, and $D_4$-invariance of the Euclidean norm yields \eqref{eq:n3-norm}.
\end{proof}

\Cref{prop:reciprocity} makes the mechanics of spoke normalization transparent: forcing unequal spokes to equal unit length ($L_{ij}=1$ or $\ell=1$) divides each coordinate vector by its irrational radius $R_{ij}(\lambda)$ or $r_{ij}$. At order $n=3$, the corner ring sits at an irrational ratio $r_2/r_1=\sqrt{2}$, so equalizing the spoke lengths breaks the cancellation between $\mathbf d_1$ and $\mathbf d_2=-\mathbf d_1$. We now examine how this irrationality mechanism operates across the three concentric rings of order $n=4$.

\subsection{Planar Balance Wheels: The \texorpdfstring{$\sqrt{5}$}{sqrt(5)} Irrationality Barrier and Power-Law Restoration}

\begin{figure}[htbp]
\centering
\includegraphics[width=\textwidth]{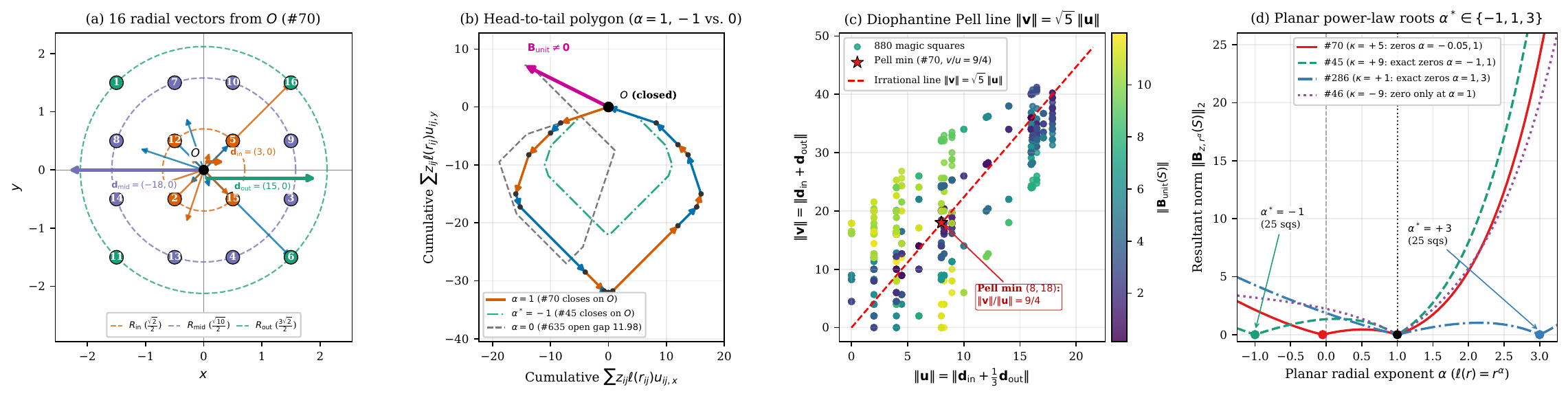}
\caption{Planar radial vectors from the origin $O=(0,0)$, head-to-tail polygon closure, and the Diophantine Pell minimum at order four. (a)~Radial vectors $z_{ij}\mathbf r_{ij}$ on square \#70 across the three concentric rings $R_{\mathrm{in}},R_{\mathrm{mid}},R_{\mathrm{out}}$, with ring dipoles $\mathbf d_{\mathrm{in}}=(3,0)$, $\mathbf d_{\mathrm{out}}=(15,0)$, and $\mathbf d_{\mathrm{mid}}=(-18,0)$. (b)~Ordered head-to-tail placement closes the $16$-gon at linear lengths $\alpha=1$ (\#70) and inverse lengths $\alpha^*=-1$ (\#45), whereas unit lengths $\alpha=0$ (\#635) leave a nonzero gap $\mathbf B_{\mathrm{unit}}(S)$ (magenta). (c)~Census of $\|\mathbf v(S)\|_2$ versus $\|\mathbf u(S)\|_2$ across all $880$ squares against the irrational balance line $\|\mathbf v\|=\sqrt{5}\,\|\mathbf u\|$, with the eight Pell-minimum squares ($v/u=9/4$) marked by red stars. (d)~Resultant norm $\|\mathbf B_{z,r^\alpha}(S)\|_2$ versus exponent $\alpha$, showing the universal zero at $\alpha=1$ and second roots at $\alpha^*=-1$ (\#45), $-0.050$ (\#70), and $+3$ (\#286).}
\label{fig:ring-vectors}
\end{figure}

At order $n=4$, the sixteen cells of $A_4=\{-\frac32,-\frac12,\frac12,\frac32\}^2$ partition into three concentric rings around the origin $O$ (\Cref{fig:hero-squares}(b), \Cref{fig:ring-vectors}(a)): four inner cells $R_{\mathrm{in}}=\{(\pm\tfrac12,\pm\tfrac12)\}$ of radius $r_1=\sqrt{2}/2$, eight middle edge-interior cells $R_{\mathrm{mid}}=\{(\pm\tfrac32,\pm\tfrac12),(\pm\tfrac12,\pm\tfrac32)\}$ of radius $r_2=\sqrt{10}/2$, and four outer corner cells $R_{\mathrm{out}}=\{(\pm\tfrac32,\pm\tfrac32)\}$ of radius $r_3=3\sqrt{2}/2$. Whereas the outer and inner radii have a rational ratio $r_3/r_1=3$, the middle and inner radii have the irrational ratio $r_2/r_1=\sqrt{5}\notin\Q$. Consequently, when spokes are trimmed to unit length ($\alpha=0$), the inner and outer rings combine into an integer lattice vector $\mathbf u(S)$ that attempts to cancel an irrational multiple $\mathbf v(S)/\sqrt{5}$ from the middle ring.

\begin{theorem}[Three-Ring Decomposition, Irrationality No-Go, and Pell Minimum]
\label{thm:three-ring}
For any arrangement $S$ of order $4$, define the three ring dipole vectors radiating from $O$ by $\mathbf d_{\mathrm{in}}(S)=\sum_{c\in R_{\mathrm{in}}}z_c\mathbf r_c$, $\mathbf d_{\mathrm{mid}}(S)=\sum_{c\in R_{\mathrm{mid}}}z_c\mathbf r_c$, and $\mathbf d_{\mathrm{out}}(S)=\sum_{c\in R_{\mathrm{out}}}z_c\mathbf r_c$. For every $4\times4$ arrangement $S$ of distinct integers, $\mathbf d_{\mathrm{in}}(S)\ne\mathbf 0$ and $\mathbf d_{\mathrm{out}}(S)\ne\mathbf 0$; for every order-four magic square $S\in\mathcal M_4$, both $R_{\mathrm{in}}$ and $R_{\mathrm{out}}$ sum to $\mu_4=34$, so $\mathbf d_{\mathrm{in}}(S)\in\Z^2\setminus\{\mathbf 0\}$ and $\mathbf d_{\mathrm{out}}(S)\in(3\Z)^2\setminus\{\mathbf 0\}$ satisfy
\begin{equation}
\label{eq:ring-lower}
\|\mathbf d_{\mathrm{in}}(S)\|_2 \ge 1,\qquad \|\mathbf d_{\mathrm{out}}(S)\|_2 \ge 3.
\end{equation}
Whenever $\sum_{c\in A_4}z_c\mathbf r_c=\mathbf 0$ (in particular, for every semi-magic square), $\mathbf d_{\mathrm{mid}}(S)=-\bigl(\mathbf d_{\mathrm{in}}(S)+\mathbf d_{\mathrm{out}}(S)\bigr)$, and for every radial length profile $\ell=(\ell_1,\ell_2,\ell_3)$ on $(r_1,r_2,r_3)$,
\begin{equation}
\label{eq:three-ring-general}
\mathbf B_{z,\ell}(S) = \left(\frac{\ell_1}{r_1} - \frac{\ell_2}{r_2}\right)\mathbf d_{\mathrm{in}}(S)
+ \left(\frac{\ell_3}{r_3} - \frac{\ell_2}{r_2}\right)\mathbf d_{\mathrm{out}}(S).
\end{equation}
In particular, for planar unit vectors ($\ell_1=\ell_2=\ell_3=1$, i.e.\ $\alpha=0$),
\begin{equation}
\label{eq:unit-uv}
\mathbf B_{\mathrm{unit}}(S) = \sqrt{2}\,\left(\mathbf u(S) - \frac{1}{\sqrt{5}}\,\mathbf v(S)\right),
\end{equation}
where $\mathbf u(S)=\mathbf d_{\mathrm{in}}(S)+\frac13\mathbf d_{\mathrm{out}}(S)\in\Z^2$ and $\mathbf v(S)=\mathbf d_{\mathrm{in}}(S)+\mathbf d_{\mathrm{out}}(S)=-\mathbf d_{\mathrm{mid}}(S)\in\Z^2$. Because $\sqrt{5}\notin\Q$ and $\mathbf d_{\mathrm{in}}(S),\mathbf d_{\mathrm{out}}(S)\ne\mathbf 0$, no $4\times4$ arrangement of distinct integers with vanishing first moments (and hence no $4\times4$ semi-magic or magic square) ever satisfies $\mathbf B_{\mathrm{unit}}(S)=\mathbf 0$, nor can $\mathbf B_{z,r^\alpha}(S)=\mathbf 0$ hold for any even integer $\alpha\in2\Z$. Across all $880$ magic squares of order~$4$, the global minimum of $\|\mathbf B_{\mathrm{unit}}(S)\|_2$ is
\begin{equation}
\label{eq:pell-min}
\min_{S\in\mathcal M_4} \|\mathbf B_{\mathrm{unit}}(S)\|_2
= \frac{18\sqrt{10}}{5} - 8\sqrt{2}
= \frac{2\sqrt{10}}{5(9+4\sqrt{5})}
= \frac{2\sqrt{10}}{5\varphi^6} \approx 0.07049108,
\end{equation}
where $\varphi=(1+\sqrt{5})/2$. This minimum is attained when $\mathbf d_{\mathrm{out}}(S)=5\mathbf d_{\mathrm{in}}(S)$ with $\|\mathbf d_{\mathrm{in}}(S)\|_2=3$, where $\mathbf u(S)$ and $\mathbf v(S)$ are parallel in the ratio $\|\mathbf v\|/\|\mathbf u\|=18/8=9/4$, the fundamental Pell convergent to $\sqrt{5}$ ($9^2-5\cdot4^2=1$). Up to $D_4$-symmetry, exactly eight magic squares achieve \eqref{eq:pell-min}: squares \#70, \#76, \#121, \#188, \#452, \#459, \#490, and \#574, all in Dudeney's Group~VI.
\end{theorem}

\begin{proof}
Labelling the four entries of $R_{\mathrm{in}}$ at $(-\frac12,\frac12)$, $(\frac12,\frac12)$, $(-\frac12,-\frac12)$, $(\frac12,-\frac12)$ by $a,b,c,d\in\Z$ gives $\mathbf d_{\mathrm{in}}(S)=\frac12\bigl((b+d)-(a+c),(a+b)-(c+d)\bigr)$. If $\mathbf d_{\mathrm{in}}(S)=\mathbf 0$, then $b+d=a+c$ and $a+b=c+d$; adding and subtracting these equations forces $b=c$ and $a=d$, contradicting that $S$ has distinct entries. For $S\in\mathcal M_4$, subtracting the sum of all sixteen cells ($4\mu_4$) from the sum of the two outer rows and two outer columns ($R_0+R_3+C_0+C_3=4\mu_4$) gives $\sum_{c\in R_{\mathrm{out}}}s_c-\sum_{c\in R_{\mathrm{in}}}s_c=0$, while the two main diagonals sum to $D_1+D_2=\sum_{c\in R_{\mathrm{out}}}s_c+\sum_{c\in R_{\mathrm{in}}}s_c=2\mu_4$; hence $\sum_{c\in R_{\mathrm{in}}}s_c=\sum_{c\in R_{\mathrm{out}}}s_c=\mu_4=34$ (an even integer), which implies $\mathbf d_{\mathrm{in}}(S)\in\Z^2\setminus\{\mathbf 0\}$. Since the corner coordinates are $(\pm\frac32,\pm\frac32)=3(\pm\frac12,\pm\frac12)$, identical algebra gives $\mathbf d_{\mathrm{out}}(S)\in(3\Z)^2\setminus\{\mathbf 0\}$ and \eqref{eq:ring-lower}. Partitioning $A_4=R_{\mathrm{in}}\sqcup R_{\mathrm{mid}}\sqcup R_{\mathrm{out}}$ gives $\mathbf d_{\mathrm{mid}}=-(\mathbf d_{\mathrm{in}}+\mathbf d_{\mathrm{out}})$, and substituting this relation into \eqref{eq:B-planar} yields \eqref{eq:three-ring-general}.

Setting $\ell_1=\ell_2=\ell_3=1$ with $1/r_1=\sqrt{2}$, $1/r_2=\sqrt{2}/\sqrt{5}$, and $1/r_3=\sqrt{2}/3$ in \eqref{eq:three-ring-general} gives $\mathbf B_{\mathrm{unit}}(S)=(\sqrt{2}-\sqrt{2}/\sqrt{5})\mathbf d_{\mathrm{in}}+(\sqrt{2}/3-\sqrt{2}/\sqrt{5})\mathbf d_{\mathrm{out}}=\sqrt{2}\bigl((\mathbf d_{\mathrm{in}}+\frac13\mathbf d_{\mathrm{out}})-\frac{1}{\sqrt{5}}(\mathbf d_{\mathrm{in}}+\mathbf d_{\mathrm{out}})\bigr)$, which is \eqref{eq:unit-uv}. Because $\mathbf u(S),\mathbf v(S)\in\Z^2$ and $\sqrt{5}\notin\Q$, $\mathbf B_{\mathrm{unit}}(S)=\mathbf 0$ requires $\mathbf u(S)=\mathbf v(S)=\mathbf 0$, which forces $\mathbf d_{\mathrm{in}}(S)=\mathbf d_{\mathrm{out}}(S)=\mathbf 0$, a contradiction. The same irrationality holds for any even exponent $\alpha=2m$, where $r_2^{\alpha-1}/r_1^{\alpha-1}=5^{m-1}\sqrt{5}\notin\Q$.

To establish the Pell minimum \eqref{eq:pell-min}, subtracting $C_1+C_2=2\mu_4$ and $C_0+C_3=2\mu_4$ from $D_1+D_2=2\mu_4$ gives $s_{00}+s_{03}=s_{31}+s_{32}$ and $s_{11}+s_{12}=s_{20}+s_{23}$; across all $24$ modulo-$2$ parity patterns of $\mathcal M_4$, the four top diagonal cells $T_{\mathrm{diag}}=\{(0,0),(0,3),(1,1),(1,2)\}$ (and likewise $R_{\mathrm{diag}}=\{(0,3),(3,3),(1,2),(2,2)\}$) always have an even sum $\sum_{c\in T_{\mathrm{diag}}}s_c\in\{18,20,\dots,42\}$. Since $u_y(S)=\sum_{c\in T_{\mathrm{diag}}}s_c-34$, $u_x(S)=\sum_{c\in R_{\mathrm{diag}}}s_c-34$, and $\mathbf v(S)-\mathbf u(S)=\frac23\mathbf d_{\mathrm{out}}(S)\in(2\Z)^2$, every magic square $S\in\mathcal M_4$ satisfies $\mathbf u(S)=2\tilde{\mathbf u}(S)\in(2\Z)^2$ and $\mathbf v(S)=2\tilde{\mathbf v}(S)\in(2\Z)^2$ with $\tilde{\mathbf u}(S),\tilde{\mathbf v}(S)\in\Z^2\setminus\{\mathbf 0\}$, $\|\tilde{\mathbf u}(S)\|_\infty\le 8$, and $\|\tilde{\mathbf u}(S)\|_2\le 4\sqrt{5}<9$. Consequently, \eqref{eq:unit-uv} becomes $\|\mathbf B_{\mathrm{unit}}(S)\|_2=2\sqrt{2}\,\|\tilde{\mathbf u}(S)-\tilde{\mathbf v}(S)/\sqrt{5}\|_2$. Decomposing $\tilde{\mathbf v}$ parallel and orthogonal to $\tilde{\mathbf u}$ shows that this norm is minimized when $\tilde{\mathbf u}=\tilde u\mathbf e$ and $\tilde{\mathbf v}=\tilde v\mathbf e$ lie along a coordinate axis $\mathbf e\in\{\pm\hat{\mathbf e}_x,\pm\hat{\mathbf e}_y\}$ ($\|\mathbf e\|_2=1$) with $1\le\tilde u\le 8$. By Lagrange's theorem on continued fractions for $\sqrt{5}=[2;\overline{4}]$ (whose principal convergents are $2/1,9/4,38/17,\dots$), across all integer denominators $1\le\tilde u\le 8<17$ the unique best rational approximation to $\sqrt{5}$ is the fundamental Pell convergent $\tilde v/\tilde u=9/4$ ($9^2-5\cdot 4^2=1$), corresponding to $(u,v)=(8,18)$ ($\mathbf d_{\mathrm{out}}=5\mathbf d_{\mathrm{in}}$ with $\|\mathbf d_{\mathrm{in}}\|_2=3$) and giving $2\sqrt{2}(9/\sqrt{5}-4)=18\sqrt{10}/5-8\sqrt{2}=2\sqrt{10}/(5\varphi^6)\approx0.07049108$. On square \#70 (\Cref{fig:hero-squares}(b)), the inner ring $\{5,12,2,15\}$ gives $\mathbf d_{\mathrm{in}}=(3,0)$ and the outer corners $\{1,16,11,6\}$ give $\mathbf d_{\mathrm{out}}=(15,0)=5\mathbf d_{\mathrm{in}}$, realizing $\mathbf u=(8,0)$, $\mathbf v=(18,0)$, and \eqref{eq:pell-min}.
\end{proof}

How, then, can a planar balance wheel achieve exact equilibrium at a non-linear length profile $\ell(r)=r^\alpha$ ($\alpha\ne1$)? Equation~\eqref{eq:three-ring-general} shows that balance is restored whenever the inner and outer ring dipoles align in the same direction ($\mathbf d_{\mathrm{out}}=\kappa\mathbf d_{\mathrm{in}}$ with $\kappa>0$), so that both rings jointly oppose the middle ring $\mathbf d_{\mathrm{mid}}=-(1+\kappa)\mathbf d_{\mathrm{in}}$ (\Cref{fig:ring-vectors}(d)).

\begin{corollary}[Restoring planar balance via power-law exponents $\alpha^*\ne1$]
\label{cor:planar-power-roots}
For power-law radial lengths $\ell(r)=r^\alpha$ ($\alpha\in\R$), equation~\eqref{eq:three-ring-general} takes the form
\begin{equation}
\label{eq:planar-power-form}
\mathbf B_{z,r^\alpha}(S) = r_1^{\alpha-1}\Bigl[\bigl(1 - 5^{(\alpha-1)/2}\bigr)\mathbf d_{\mathrm{in}}(S) + \bigl(3^{\alpha-1} - 5^{(\alpha-1)/2}\bigr)\mathbf d_{\mathrm{out}}(S)\Bigr].
\end{equation}
Because $1<\sqrt{5}<3$, the ratio $h(\alpha)=\frac{5^{(\alpha-1)/2}-1}{3^{\alpha-1}-5^{(\alpha-1)/2}}$ (extended by $h(1)=\frac{\ln 5}{\ln(9/5)}\approx2.7381$) is a strictly decreasing bijection from $\R$ onto $(0,\infty)$. Consequently, an order-four magic square $S$ possesses a second power-law balance exponent $\alpha^*(S)\ne1$ satisfying $\mathbf B_{z,r^{\alpha^*}}(S)=\mathbf 0$ if and only if $\mathbf d_{\mathrm{in}}(S)$ and $\mathbf d_{\mathrm{out}}(S)$ are strictly parallel ($\mathbf d_{\mathrm{out}}(S)=\kappa\mathbf d_{\mathrm{in}}(S)$ with $\kappa>0$), which holds for $204$ of the $880$ magic squares (while $204$ squares have anti-parallel dipoles $\kappa<0$). Across all $880$ magic squares, the second root $\alpha^*(\kappa)$ is an integer in only two cases: inverse radial lengths $\alpha^*=-1$ ($\ell(r)=1/r$, corresponding to $\kappa=h(-1)=+9$), achieved by $25$ squares including \#45 (\Cref{fig:hero-squares}(c)), and cubic radial lengths $\alpha^*=+3$ ($\ell(r)=r^3$, corresponding to $\kappa=h(3)=+1$), achieved by $25$ squares including \#286. For the $17$ squares with $\kappa=+5$ (including the eight Pell-minimum squares of \Cref{thm:three-ring}), the exact zero occurs at $\alpha^*(5)\approx -0.0496$, explaining why $\|\mathbf B_{z,r^\alpha}(S)\|_2$ pinches nearly to zero at $\alpha=0$ (\Cref{fig:ring-vectors}(d)).
\end{corollary}

\begin{proof}
Setting \eqref{eq:planar-power-form} to $\mathbf 0$ for $\alpha\ne1$ requires $\mathbf d_{\mathrm{out}}(S)=h(\alpha)\mathbf d_{\mathrm{in}}(S)$. Since $1<5^{t/2}<3^t$ for $t=\alpha-1>0$ and $3^t<5^{t/2}<1$ for $t<0$, both numerator and denominator of $h(\alpha)$ have the sign of $\alpha-1$, so $h(\alpha)>0$ with $\lim_{\alpha\to-\infty}h(\alpha)=+\infty$ and $\lim_{\alpha\to+\infty}h(\alpha)=0$. Since $\kappa\in\Q_{>0}$ cannot equal the transcendental number $\frac{\ln 5}{\ln(9/5)}$, every square with $\kappa>0$ has a unique root $\alpha^*(\kappa)=h^{-1}(\kappa)\ne1$. By \Cref{thm:three-ring}, an integer root $\alpha^*\ne1$ must be odd, $\alpha^*=2m+1$ ($m\in\Z\setminus\{0\}$), where $h(2m+1)=\frac{5^m-1}{9^m-5^m}$; comparing this sequence against the $15$ positive rational ratios $\kappa$ realized on $\mathcal M_4$ leaves only $m=-1$ ($\kappa=9$, $25$ squares) and $m=+1$ ($\kappa=1$, $25$ squares).
\end{proof}

\subsection{Three-Dimensional Spoke Stars: Unit-Sphere Counterweights and Cubic Vertical Balance}

\begin{figure}[htbp]
\centering
\includegraphics[width=\textwidth]{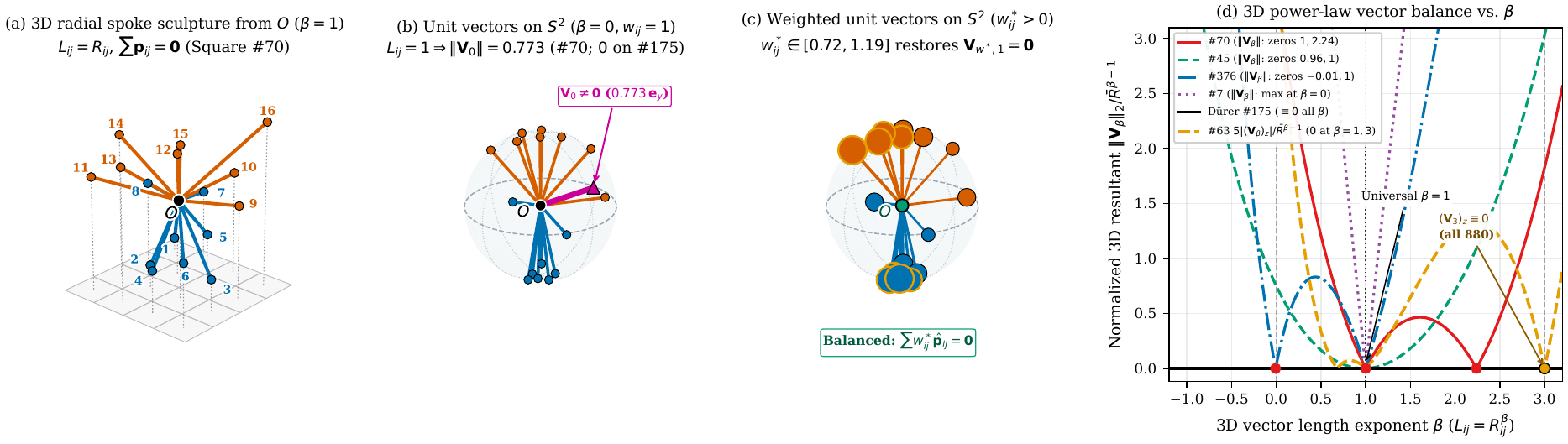}
\caption{Three-dimensional radial spoke sculptures from the origin $O=(0,0,0)$ under variable length and tip weight. (a)~Natural $16$-spoke sculpture $\mathbf p_{ij}$ ($\beta=1$, orange $z_{ij}>0$, blue $z_{ij}<0$) of non-associative square \#70, balancing identically at $O$. (b)~Trimming all spokes to the unit sphere $S^2$ ($\beta=0$, equal weights) leaves a non-zero resultant $\mathbf V_0=(0,0.773,0)$ (magenta arrow) on \#70, whereas antipodal pairs cancel on D\"urer \#175. (c)~Minimum-variance positive tip weights $w_{ij}^*\in[0.72,1.19]$ (gold-rimmed for $w_{ij}^*>1$) restore exact unit-sphere balance $\sum w_{ij}^*\hat{\mathbf p}_{ij}=\mathbf 0$. (d)~Normalized resultant $\|\mathbf V_\beta(S;1)\|_2/\bar R^{\beta-1}$ versus exponent $\beta$, showing universal balance at $\beta=1$, second roots $\beta^*\ne1$ on Group~VI squares (\#70, \#45, \#376), and universal vertical vanishing $(\mathbf V_3)_z=0$ at $\beta=3$ (gold dashed curve).}
\label{fig:origin-vectors-3d}
\end{figure}

We now lift the spoke sculpture into three dimensions (\Cref{fig:origin-vectors-3d}). In 3D, trimming all sixteen spokes to the unit sphere $S^2$ ($\beta=0$) breaks equilibrium on $832$ of the $880$ squares, sparing only the $48$ associative squares of Group~III where antipodal pairs cancel pairwise ($\hat{\mathbf p}_{-c}=-\hat{\mathbf p}_c$). Nevertheless, because the natural unequal lengths $R_{ij}(\lambda)>0$ satisfy $\sum_{i,j}R_{ij}(\lambda)\hat{\mathbf p}_{ij}(\lambda)=\mathbf 0$, the origin $O$ always lies strictly inside the convex hull of the sixteen unit rays, guaranteeing that positive unequal tip weights $\mathbf w^*>\mathbf 0$ restore unit-sphere equilibrium for every magic square.

\begin{theorem}[3D Power-Law Vectors, Cubic Vertical Vanishing, and Positive Unit-Sphere Weights]
\label{thm:3d-power}
Let $S$ be a semi-magic square of order $n$ and let $\lambda>0$. For cubic 3D spoke lengths $L_{ij}=R_{ij}(\lambda)^3$ from the origin (dual by \Cref{prop:reciprocity} to weighting $\mathbf p_{ij}(\lambda)$ by $w_{ij}=R_{ij}(\lambda)^2=r_{ij}^2+\lambda^2z_{ij}^2$), the vertical $z$-component vanishes identically for every semi-magic square of any order $n$ and every $\lambda>0$:
\begin{equation}
\label{eq:V3-exact}
\mathbf V_3(S;\lambda) = \sum_{i,j}\bigl(r_{ij}^2+\lambda^2z_{ij}^2\bigr)(x_j,y_i,\lambda z_{ij})
= \biggl(\lambda^2\sum_{i,j}z_{ij}^2 x_j,\;\; \lambda^2\sum_{i,j}z_{ij}^2 y_i,\;\; 0\biggr),
\end{equation}
and at order $n=4$, $\mathbf V_3(S;\lambda)=\mathbf 0$ if and only if $S$ is one of the $48$ associative magic squares. Similarly, under any even weight function $w_{ij}=f(|z_{ij}|)$, every associative square satisfies $\hat{\mathbf p}_{-c}(\lambda)=-\hat{\mathbf p}_c(\lambda)$ and hence $\mathbf V_\beta(S;\lambda)=\mathbf 0$ for all $\beta\in\R$ and $\lambda>0$; at $n=4$, $\mathbf V_0(S;\lambda)=\mathbf 0$ holds if and only if $S$ is one of the $48$ associative squares ($\|\mathbf V_0(S;1)\|_2\in[0.0072,2.4284]$ on the remaining $832$ squares), while among the $304$ Group~VI squares (where $C_2$ symmetry confines $\mathbf V_\beta(S;\lambda)$ to a single coordinate axis), $107$ squares possess a second 3D power-law balance exponent $\beta^*(S)\in[-3,4]\setminus\{1\}$ with $\mathbf V_{\beta^*}(S;1)=\mathbf 0$ (\Cref{fig:origin-vectors-3d}(d)). Finally, letting $\hat P(\lambda)\in\R^{N\times 3}$ have rows $\hat{\mathbf p}_{ij}(\lambda)^{\top}$, the unique weight vector $\mathbf w^*(S;\lambda)$ of minimum Euclidean distance $\|\mathbf w-\mathbf 1\|_2$ from uniform weights $\mathbf 1$ satisfying $\sum_{i,j}w_{ij}^*\hat{\mathbf p}_{ij}(\lambda)=\mathbf 0$ is
\begin{equation}
\label{eq:w-star}
\mathbf w^*(S;\lambda) = \mathbf 1 - \hat P(\lambda)\bigl(\hat P(\lambda)^{\top}\hat P(\lambda)\bigr)^{-1}\mathbf V_0(S;\lambda),
\end{equation}
which is strictly positive ($\mathbf w^*(S;\lambda)>\mathbf 0$) across all $880$ order-four magic squares at $\lambda=1$, $\lambda_{\mathrm{iso}}=1/\sqrt{17}$, and $\lambda_{\mathrm{cube}}=0.20$ (\Cref{fig:origin-vectors-3d}(c)), with $\|\mathbf w^*(S;1)-\mathbf 1\|_2\in[0.0052,\,1.5967]$ across the $832$ non-associative squares.
\end{theorem}

\begin{proof}
Expanding $\mathbf V_3(S;\lambda)=\sum_{i,j}(r_{ij}^2+\lambda^2z_{ij}^2)(x_j,y_i,\lambda z_{ij})$, the horizontal sum $\sum_{i,j}r_{ij}^2(x_j,y_i)$ vanishes by odd grid symmetry, leaving $\lambda^2\sum_{i,j}z_{ij}^2(x_j,y_i)$. The vertical component is $\lambda\sum_{i,j}(x_j^2+y_i^2)z_{ij}+\lambda^3\sum_{i,j}z_{ij}^3$; here $\sum_{i,j}(x_j^2+y_i^2)z_{ij}=0$ because $x_j^2+y_i^2\in W_{\mathrm{rc}}$, and $\sum_{i,j}z_{ij}^3=0$ by odd symmetry of $\{z_{ij}\}$, proving \eqref{eq:V3-exact}. For associative squares, $z(-x,-y)=-z(x,y)$ implies $\mathbf p_{-c}(\lambda)=-\mathbf p_c(\lambda)$ and $R_{-c}(\lambda)=R_c(\lambda)$, so opposite cells cancel pairwise for every $\beta\in\R$. Minimizing $\frac12\|\mathbf w-\mathbf 1\|_2^2$ subject to $\hat P(\lambda)^{\top}\mathbf w=\mathbf 0$ by Lagrange multipliers gives \eqref{eq:w-star}, and exact enumeration over the $880$ squares establishes positivity and the stated bounds.
\end{proof}

Notice how the cubic resultant $\mathbf V_3(S;\lambda)$ in \eqref{eq:V3-exact} bridges static rays and rotational dynamics: while its vertical component vanishes for every semi-magic square, its residual horizontal dipole $\lambda^2\sum_{i,j}z_{ij}^2(x_j,y_i)$ measures the quadratic height asymmetry across the grid. In \Cref{sec:dynamics}, we set our sculptures into three-dimensional rotation and show that this exact vector governs the vertical off-diagonal products of inertia of an entry-weighted gem.

\section{Rotational Dynamics: Smooth Precession and Dzhanibekov Flips}
\label{sec:dynamics}

Our second test moves from static radial equilibrium to three-dimensional rotation. Suppose a Magic Gem or radial spoke sculpture of total mass $M$ and centre of mass $\mathbf c$ is supported on a low-friction conical pivot or gimbal at the origin $O$ under uniform gravity $\mathbf g=(0,0,-g)^{\top}$ (generalizing the conical fulcrum of \texttt{weighted\_ring.stl} in \cite{mathewson2025magicgemweb}), or tossed in free space. About the pivot $O$, gravity exerts the torque
\[
\boldsymbol\tau_{\mathrm{grav}} = M\,\mathbf c \times \mathbf g,
\]
which vanishes identically in every 3D orientation if and only if the centre of mass lies at the pivot ($\mathbf c=\mathbf 0$). Consequently, any balanced sculpture ($\mathbf c=\mathbf 0$) supported at $O$ rotates under Earth's gravity according to Euler's exact torque-free equations \cite{arnold1989mathematical,goldstein2002classical}: in the principal body frame, the angular velocity vector $\boldsymbol\omega=(\omega_1,\omega_2,\omega_3)^{\top}$ evolves via
\begin{equation}
\label{eq:euler}
I_1\dot\omega_1 = (I_2-I_3)\omega_2\omega_3,\qquad
I_2\dot\omega_2 = (I_3-I_1)\omega_3\omega_1,\qquad
I_3\dot\omega_3 = (I_1-I_2)\omega_1\omega_2,
\end{equation}
where $I_1,I_2,I_3$ are the principal moments of inertia. For an asymmetric body with three distinct principal moments $I_1<I_2<I_3$, rotation about the intermediate axis $I_2$ is exponentially unstable: any small perturbation grows at rate $\sqrt{(I_2-I_1)(I_3-I_2)/(I_1I_3)}\,|\omega_2|$, causing the body to execute the $180^\circ$ flips of the Dzhanibekov (or tennis-racket) theorem \cite{ashbaugh1991twisting,poinsot1834theorie}.

When the sixteen unequal heights $z_{ij}$ of a magic square are set spinning in $\R^3$, do their geometric inequalities produce an asymmetric top that tumbles via Dzhanibekov flips, or does line-sum equality enforce smooth axisymmetric precession?

\subsection{Equal Masses at Unequal Heights Form a Universal Axisymmetric Rotor}

Consider first the discrete point cloud $P_S(\lambda)=\{\mathbf p_{ij}(\lambda)=(x_j,y_i,\lambda z_{ij})\}$ carrying equal point masses $m=M/N$ at the $N=n^2$ unequal lifted points (for instance, identical steel bearings mounted at the tips of the sixteen radial spokes $\mathbf p_{ij}(\lambda)$ on a lightweight hub at $O$). Because $\mathbf c=\frac1N\sum_{i,j}\mathbf p_{ij}(\lambda)=\mathbf 0$, gravitational torque about $O$ vanishes identically for every semi-magic square, and the magic line sums neutralize all off-diagonal products of inertia.

\begin{figure}[htbp]
\centering
\includegraphics[width=\textwidth]{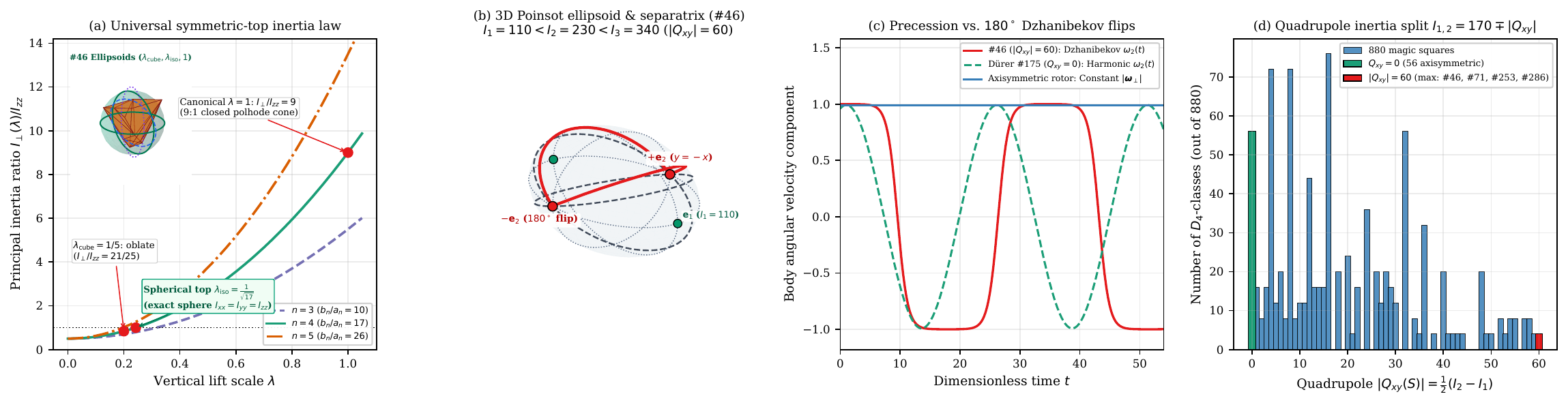}
\caption{Three-dimensional rotational dynamics of Magic Gems. (a)~Universal inertia ratio $I_\perp(\lambda)/I_{zz}=\frac12(1+\lambda^2(n^2+1))$ versus vertical scale $\lambda$ for $n=3,4,5$, with 3D inertia ellipsoids at $\lambda_{\mathrm{cube}}=1/5$ ($21/25$), the spherical-top transition $\lambda_{\mathrm{iso}}=1/\sqrt{17}$ ($1$), and $\lambda=1$ ($9$) for $n=4$. (b)~Triaxial Poinsot energy ellipsoid of entry-weighted square \#46 ($I_1=110<I_2=230<I_3=340$), showing the figure-eight separatrix and heteroclinic polhode connecting $+\omega_2$ to $-\omega_2$. (c)~Euler time series \eqref{eq:euler}: an axisymmetric magic rotor precesses harmonically with constant $|\boldsymbol\omega_\perp|$, whereas entry-weighted square \#46 ($|Q_{xy}|=60$) executes periodic $180^\circ$ Dzhanibekov flips between $+\omega_2$ and $-\omega_2$. (d)~Census of the quadrupole split $|Q_{xy}(S)|=\frac12(I_2-I_1)$ across all $880$ order-four magic squares.}
\label{fig:precession}
\end{figure}

\begin{theorem}[Universal symmetric-top law, spherical-top scale, and periodic polhodes]
\label{thm:rotor}
Let $S$ be any arrangement of $\{1,\dots,n^2\}$ ($n\ge2$) whose height field has vanishing first moments $\sum_{i,j}x_jz_{ij}=\sum_{i,j}y_iz_{ij}=0$ (in particular, any semi-magic or magic square), and place equal point masses $m=M/N$ at the points $\mathbf p_{ij}(\lambda)=(x_j,y_i,\lambda z_{ij})$. For every $\lambda>0$, the inertia tensor about the origin $O$,
\[
I(\lambda) = \frac{M}{N}\sum_{i,j}\Bigl(\|\mathbf p_{ij}(\lambda)\|_2^2\,I_3 - \mathbf p_{ij}(\lambda)\mathbf p_{ij}(\lambda)^{\top}\Bigr),
\]
is diagonal and strictly axisymmetric about the $z$-axis:
\begin{equation}
\label{eq:inertia-point}
I(\lambda) = M\diag\bigl(I_\perp(\lambda),\, I_\perp(\lambda),\, I_{zz}\bigr),\qquad
I_\perp(\lambda) = a_n + \lambda^2 b_n,\qquad
I_{zz} = 2a_n,
\end{equation}
with universal principal inertia ratio
\begin{equation}
\label{eq:inertia-ratio}
\frac{I_\perp(\lambda)}{I_{zz}} = \frac{a_n + \lambda^2 b_n}{2a_n} = \frac{1 + \lambda^2(n^2+1)}{2}.
\end{equation}
Varying the vertical lift scale $\lambda$ drives the rotor through three critical aspect ratios (\Cref{fig:precession}(a)). At the canonical lift $\lambda=1$,
\begin{equation}
\label{eq:ratio-lambda1}
\frac{I_\perp(1)}{I_{zz}} = \frac{n^2+2}{2},
\end{equation}
which is an exact integer for every even order $n$ ($9$ for $n=4$, $19$ for $n=6$, $33$ for $n=8$) and a half-integer for every odd order $n$ ($11/2$ for $n=3$, $27/2$ for $n=5$). At the spherical-top critical scale
\begin{equation}
\label{eq:lambda-iso}
\lambda_{\mathrm{iso}}(n) = \frac{1}{\sqrt{n^2+1}}
\qquad\Bigl(\tfrac{1}{\sqrt{10}}\approx0.3162\text{ for }n=3,\quad \tfrac{1}{\sqrt{17}}\approx0.2425\text{ for }n=4\Bigr),
\end{equation}
$I_\perp(\lambda_{\mathrm{iso}})=I_{zz}=2a_n$, so $I(\lambda_{\mathrm{iso}})=2Ma_nI_3$ is a scalar multiple of the identity and the point-weighted gem is a \emph{spherical top} for which every axis through $O$ in $\R^3$ is a principal axis of inertia. At the cubic bounding-box scale $\lambda_{\mathrm{cube}}(n)=\frac{n-1}{n^2-1}=\frac{1}{n+1}$ (where the vertical span matches the horizontal footprint $[-\frac{n-1}{2},\frac{n-1}{2}]^2$, giving $\lambda_{\mathrm{cube}}=1/5=0.20$ at $n=4$),
\begin{equation}
\label{eq:ratio-cube}
\frac{I_\perp(\lambda_{\mathrm{cube}})}{I_{zz}} = \frac{n^2+n+1}{(n+1)^2}
\qquad\Bigl(\frac{21}{25}=0.84\text{ at }n=4\Bigr).
\end{equation}
Because $I_{xx}(\lambda)=I_{yy}(\lambda)=I_\perp(\lambda)$, no distinct intermediate principal axis exists: in torque-free motion \eqref{eq:euler}, $\omega_3(t)=\omega_{3,0}$ is constant and the transverse angular velocity $(\omega_1(t),\omega_2(t))$ precesses in a uniform circle at the body-frame polhode frequency
\begin{equation}
\label{eq:polhode-freq}
\Omega_{\mathrm{body}}(\lambda) = \frac{I_\perp(\lambda)-I_{zz}}{I_\perp(\lambda)}\,\omega_{3,0}
= \frac{\lambda^2(n^2+1)-1}{\lambda^2(n^2+1)+1}\,\omega_{3,0},
\end{equation}
closing periodically after $(n^2+2)/2=9$ spin revolutions at $\lambda=1$ ($n=4$) and rotating without wobble ($\Omega_{\mathrm{body}}=0$) about any initial angular velocity vector $\boldsymbol\omega_0\in\R^3$ at $\lambda=\lambda_{\mathrm{iso}}(n)$.
\end{theorem}

\begin{proof}
By \cite[Lem.~3.1]{mathewson2026gems}, the empirical covariance matrix $\Sigma(\lambda)=\frac1N\sum_{i,j}\mathbf p_{ij}(\lambda)\mathbf p_{ij}(\lambda)^{\top}$ has diagonal entries $\Sigma_{xx}=\Sigma_{yy}=a_n=(n^2-1)/12$ and $\Sigma_{zz}=\lambda^2b_n=\lambda^2(n^4-1)/12=\lambda^2(n^2+1)a_n$, while $\Sigma_{xy}=0$ by grid symmetry and $\Sigma_{xz}=\Sigma_{yz}=0$ by vanishing of the first moments. Since $I(\lambda)=M\bigl(\tr(\Sigma(\lambda))I_3-\Sigma(\lambda)\bigr)$, we obtain \eqref{eq:inertia-point} and \eqref{eq:inertia-ratio}. Evaluating \eqref{eq:inertia-ratio} at $\lambda=1$, $\lambda=(n^2+1)^{-1/2}$, and $\lambda=(n+1)^{-1}$ yields \eqref{eq:ratio-lambda1}--\eqref{eq:ratio-cube}. Substituting $I_1=I_2=I_\perp(\lambda)$ and $I_3=I_{zz}$ into Euler's equations \eqref{eq:euler} gives $\dot\omega_3=0$ and $\dot\omega_1=\Omega_{\mathrm{body}}\omega_2$, $\dot\omega_2=-\Omega_{\mathrm{body}}\omega_1$ with $\Omega_{\mathrm{body}}$ as in \eqref{eq:polhode-freq}.
\end{proof}

\subsection{Unequal Entry Masses Expose \texorpdfstring{$Q_{xy}$}{Qxy} and Trigger Dzhanibekov Flips}

The universal axisymmetric rotation of \Cref{thm:rotor} breaks as soon as each cell is weighted by its \emph{unequal magic entry} $m_{ij}=s_{ij}$, as in the planar balancing wheel \texttt{weighted\_ring.stl} \cite{mathewson2025magicgemweb}. While Loly's theorem \cite{loly2004invariance} keeps the diagonal moments equal ($I_{xx}=I_{yy}=I_0$), the unequal masses expose a second-order quadrant inequality—the planar quadrupole moment $Q_{xy}(S)=\sum_{i,j}z_{ij}x_jy_i$—together with the cubic 3D vector resultant $\mathbf V_3(S;\lambda)$ from \Cref{sec:vectors}. On square \#46 (\Cref{fig:hero-squares}(d)), for example, the heavy entries $\{16,13,15,14\}$ concentrate in the NW/SE quadrants ($xy<0$) while the light entries $\{1,4,2,3\}$ sit in the NE/SW quadrants ($xy>0$). Flicking the board to spin about its intermediate diagonal axis $y=-x$ ($I_2=230$) places the angular velocity directly on the saddle of the Poinsot energy ellipsoid (\Cref{fig:precession}(b)), triggering a $180^\circ$ Dzhanibekov flip (\Cref{fig:precession}(c)).

\begin{theorem}[Entry-weighted inertia split, quadrupole Dzhanibekov instability, and 3D lifted tensor]
\label{thm:quadrupole-split}
Let $S$ be any semi-magic square of order $n$, and assign entry masses $m_{ij}=s_{ij}=\bar s+z_{ij}$ with total mass $M=\sum_{i,j}s_{ij}=N\bar s=n^2(n^2+1)/2$. When the masses $m_{ij}=s_{ij}$ are placed at the planar grid positions $(x_j,y_i,0)\in\R^3$, the centre of mass lies at $O$ and the $3\times3$ inertia tensor about $O$ is
\begin{equation}
\label{eq:inertia-entry}
I^{(s)}(S) = \begin{pmatrix}
I_0 & -Q_{xy}(S) & 0\\
-Q_{xy}(S) & I_0 & 0\\
0 & 0 & 2I_0
\end{pmatrix},\quad
I_0 = M a_n = \frac{n^2(n^4-1)}{24},\quad
Q_{xy}(S) = \sum_{i,j}z_{ij}x_jy_i,
\end{equation}
whose principal moments of inertia along the horizontal diagonals $y=\pm x$ and the vertical $z$-axis are
\begin{equation}
\label{eq:principal-entry}
I_1(S) = I_0 - |Q_{xy}(S)|,\qquad
I_2(S) = I_0 + |Q_{xy}(S)|,\qquad
I_3 = 2I_0.
\end{equation}
At order $n=4$ ($I_0=170$, $I_3=340$), exactly $56$ of the $880$ magic squares satisfy $Q_{xy}(S)=0$ ($48$ associative squares in Group~III, including D\"urer's \#175, plus $8$ in Groups~IX--X), remaining axisymmetric rotors with $I_1=I_2=170$; the remaining $824$ squares have $|Q_{xy}(S)|\in\{1,\dots,60\}$ (peaking at $|Q_{xy}|=60$, $I_1=110 < I_2=230 < I_3=340$, on Group~VI squares \#46, \#71, \#253, and \#286) and undergo unstable Dzhanibekov flips about $y=-\operatorname{sgn}(Q_{xy})\,x$ at dimensionless exponential growth rate
\[
\Gamma(S) = \sqrt{\frac{(I_2-I_1)(I_3-I_2)}{I_1 I_3}} = \sqrt{\frac{|Q_{xy}(S)|}{I_0}} = \sqrt{\frac{|Q_{xy}(S)|}{170}}.
\]

When the entry masses $m_{ij}=s_{ij}=\bar s+z_{ij}$ are placed at the 3D lifted points $\mathbf p_{ij}(\lambda)=(x_j,y_i,\lambda z_{ij})$ ($\lambda>0$), larger entries $s_{ij}>\bar s$ sit at positive heights $z_{ij}>0$ while smaller entries $s_{ij}<\bar s$ sit at negative heights $z_{ij}<0$, biasing the vertical moment upward by $\lambda\sum_{i,j}z_{ij}^2$ and shifting the centre of mass to $\mathbf c^{(s)}(\lambda)=(0,0,z_{\mathrm{cm}}(\lambda))$ with $z_{\mathrm{cm}}(\lambda)=\lambda b_n/\bar s=\frac{\lambda(n^2-1)}{6}$ ($=\frac52\lambda$ at $n=4$). Because $\sum_{i,j}z_{ij}^3=0$, the diagonal moments of inertia about $O$ are $I_{xx,O}^{(s)}(\lambda)=I_{yy,O}^{(s)}(\lambda)=M(a_n+\lambda^2 b_n)$ and $I_{zz,O}^{(s)}(\lambda)=2I_0$, matching the equal-mass gem in \Cref{thm:rotor}. Translating to the true centre of mass $\mathbf c^{(s)}(\lambda)$ via Steiner's parallel-axis theorem leaves $I_{zz,\mathrm{cm}}^{(s)}=2I_0$ and all off-diagonal products unchanged while subtracting $Mz_{\mathrm{cm}}(\lambda)^2$ from the horizontal diagonal moments, giving
\[
I_{xx,\mathrm{cm}}^{(s)}(\lambda) = I_{yy,\mathrm{cm}}^{(s)}(\lambda)
= M\!\left(a_n + \lambda^2 b_n\!\left(1 - \frac{b_n}{\bar s^2}\right)\right)
= M a_n\!\left(1 + \frac{2(n^2+2)}{3}\lambda^2\right)
\]
($=170(1+12\lambda^2)$ at $n=4$). In both frames ($O$ and $\mathbf c^{(s)}(\lambda)$), the off-diagonal products of inertia are $I_{xy}^{(s)}=-Q_{xy}(S)$ and
\[
\bigl(I_{xz}^{(s)}(\lambda),\, I_{yz}^{(s)}(\lambda)\bigr)
= -\lambda\sum_{i,j}z_{ij}^2(x_j,y_i)
= -\lambda^{-1}\Bigl(\bigl(\mathbf V_3(S;\lambda)\bigr)_x,\, \bigl(\mathbf V_3(S;\lambda)\bigr)_y\Bigr),
\]
so an order-four magic square has a diagonal, axisymmetric 3D entry-weighted inertia tensor if and only if it is one of the $48$ associative squares ($Q_{xy}=0$ and $\mathbf V_3=\mathbf 0$).
\end{theorem}

\begin{proof}
Since every row excess $R_i=\sum_jz_{ij}$ and column excess $C_j=\sum_iz_{ij}$ vanishes, $\sum_{i,j}s_{ij}x_j^2=\bar s\sum_{i,j}x_j^2+\sum_jx_j^2C_j=M a_n=I_0$ and $\sum_{i,j}s_{ij}y_i^2=I_0$, giving $I_{xx}=I_{yy}=I_0$ and $I_{zz}=2I_0$. Since $\sum_{i,j}x_jy_i=0$, $I_{xy}=-\sum_{i,j}s_{ij}x_jy_i=-Q_{xy}(S)$. Diagonalizing $\bigl(\begin{smallmatrix}I_0 & -Q_{xy}\\ -Q_{xy} & I_0\end{smallmatrix}\bigr)$ yields eigenvalues $I_0\mp|Q_{xy}|$ along $(1,\pm1)/\sqrt{2}$, proving \eqref{eq:inertia-entry}--\eqref{eq:principal-entry}. At $n=4$, $I_0=170$, $|Q_{xy}|\le 60 < 170$, and $I_3-I_2=2I_0-(I_0+|Q_{xy}|)=I_1$ simplifies $\sqrt{(I_2-I_1)(I_3-I_2)/(I_1I_3)}$ to $\sqrt{|Q_{xy}|/I_0}$.

At the 3D lifted points $\mathbf p_{ij}(\lambda)$, $\sum_{i,j}s_{ij}(x_j,y_i)=\mathbf 0$ while $\sum_{i,j}(\bar s+z_{ij})\lambda z_{ij}=\lambda\sum_{i,j}z_{ij}^2=N\lambda b_n$, giving $z_{\mathrm{cm}}(\lambda)=\lambda b_n/\bar s$. Since $\sum_{i,j}z_{ij}^3=0$, $\sum_{i,j}(\bar s+z_{ij})\lambda^2z_{ij}^2=M\lambda^2b_n$, so $I_{xx,O}^{(s)}(\lambda)=I_{yy,O}^{(s)}(\lambda)=M(a_n+\lambda^2b_n)$. Steiner's theorem subtracts $Mz_{\mathrm{cm}}(\lambda)^2\diag(1,1,0)$, and $1-b_n/\bar s^2=\frac{2(n^2+2)}{3(n^2+1)}$ gives $b_n(1-b_n/\bar s^2)=\frac{2(n^2+2)}{3}a_n$. Finally, $I_{xz}^{(s)}(\lambda)=-\lambda\sum_{i,j}z_{ij}^2x_j=-\lambda^{-1}(\mathbf V_3)_x$ and similarly for $I_{yz}^{(s)}(\lambda)$, both of which vanish alongside $Q_{xy}(S)$ on all $48$ associative squares by \Cref{thm:3d-power}.
\end{proof}

\subsection{From Point Clouds to Continuous Polyhedra: Why Solid Hulls Tilt on the Pivot}

What happens to rotational stability when we replace the discrete point cloud $P_S$ by the uniform-density 3D solid polyhedron $\mathcal G_0(S)=\conv P_S$? Decomposing $\mathcal G_0(S)$ into exact rational tetrahedra (\Cref{sec:appendix}) shows that about its own centre of mass $\mathbf c_0(S)$, the solid inertia tensor remains remarkably close to the axisymmetric point-mass rotor of \Cref{thm:rotor}: across all $880$ squares at $\lambda=1$, the ratio of the mean transverse moment $\frac12(I_2+I_3)$ to the axial moment $I_1$ averages $9.41$ (ranging over $[7.34,\,11.57]$, tracking the point-mass value $9$), and the relative transverse splitting between $I_2$ and $I_3$ remains below $2.5\%$ ($\frac{I_3-I_2}{I_2}\in[0.00004,\,0.0243]$, mean $0.0071$).

Yet when we place the cast solid $\mathcal G_0(S)$ on our pivot at the grid origin $O$ under gravity, a more basic failure occurs: on $829$ of the $880$ squares, taking the convex hull shifts the solid centre of mass $\mathbf c_0(S)$ off the origin ($\mathbf c_0(S)\ne\mathbf 0$), so gravity immediately exerts a tipping torque $\boldsymbol\tau_{\mathrm{grav}}=M\mathbf c_0(S)\times\mathbf g\ne\mathbf 0$. Understanding why the convex hull breaks centroid balance—and how the interior points of the magic square restore it—is the subject of \Cref{sec:nested}.

\section{Convex Hulls and ``Diamond-in-a-Diamond'' Counterweighting}
\label{sec:nested}

Why does passing from the discrete point cloud $P_S$ ($\sum_{i,j}\mathbf p_{ij}=\mathbf 0$) to its uniform-density solid convex hull $\mathcal G_0(S):=\mathcal G(S)=\conv P_S$ shift the solid centre of mass
\[
\mathbf c_0(S) := \mathbf c_{\mathrm{solid}}(\mathcal G_0(S)) = \frac{1}{\vol\mathcal G_0(S)}\iiint_{\mathcal G_0(S)}\mathbf x\,dV
\]
away from the origin $O$ on $829$ of the $880$ order-four magic squares? Taking the three-dimensional convex hull introduces two geometric inequalities: first, the boundary $\partial\mathcal G_0(S)$ depends solely on the $V_0=|\vertset\mathcal G_0(S)|$ extreme vertices and omits the $V_1=N-V_0$ interior and non-vertex points $P_S^{(1)}=P_S\setminus\vertset\mathcal G_0(S)$; second, even when $V_0=N=16$ ($18$ of the $880$ squares), decomposing $\mathcal G_0(S)$ into tetrahedra from $O$ weights each simplex $\Delta_k$ by its volume $\vol(\Delta_k)$ rather than weighting the sixteen vertices equally.

Since the outer hull vertices and the discarded inner points together sum to zero, their centroids must pull in strictly opposite directions across $O$. We now test whether physically re-engaging the inner core $\mathcal G_1(S)=\conv P_S^{(1)}$ restores exact balance at $O$ across cast solids, 3D-printed wireframes, complete strut networks, and folded papercraft shells.

\subsection{The 51 Naturally Balanced Solids and Outer--Inner Vertex Opposition}

We first classify the $51$ order-four magic squares whose uniform solid hull $\mathcal G_0(S)$ balances at $O$ without any counterweight, and then prove the outer--inner opposition law for the remaining $829$ squares.

\begin{proposition}[Classification of the 51 naturally balanced order-four gems]
\label{prop:51-balanced}
Across the $880$ magic squares of order~$4$, the uniform solid $\mathcal G_0(S)$ has centre of mass $\mathbf c_0(S)=\mathbf 0$ for exactly $51$ squares, consisting of the $48$ associative squares (Dudeney's Group~III, where $z(-x,-y)=-z(x,y)$ implies $\mathcal G_0(S)=-\mathcal G_0(S)$ and forces $\mathbf c_0(S)=-\mathbf c_0(S)=\mathbf 0$ by central inversion) together with three non-associative squares in Dudeney's Group~VI: squares \#86, \#739, and \#756 (\Cref{tab:three-exceptions}).
\end{proposition}

\begin{proof}
For every associative square $S$, $P_S=-P_S$ implies $\mathcal G_0(S)=-\mathcal G_0(S)$ and $\mathbf c_0(S)=\mathbf 0$. Each of the three Group~VI squares \#86, \#739, and \#756 satisfies $\overline S = hS$ or $\overline S = vS$, so $\mathcal G_0(S)$ is invariant under a $180^\circ$ half-turn $\iota_h\colon(x,y,z)\mapsto(x,-y,-z)$ or $\iota_v\colon(x,y,z)\mapsto(-x,y,-z)$ in $\R^3$ \cite[Prop.~4.3]{mathewson2026gems}. The half-turn forces two coordinates of $\mathbf c_0(S)$ to vanish identically, and exact simplicial integration (\Cref{sec:appendix}) certifies that the single remaining rational axial moment vanishes exclusively on \#86, \#739, and \#756 among all $832$ non-associative squares.
\end{proof}

\begin{table}[htbp]
\centering
\caption{The three non-associative order-four magic squares whose uniform solid $\mathcal G_0(S)$ has exact centre of mass $\mathbf c_0(S)=\mathbf 0$, together with their inner core $\mathcal G_1(S)=\conv(P_S\setminus\vertset\mathcal G_0(S))$. All three squares lie in Dudeney's Group~VI and have $12$ outer vertices and $4$ inner points (spanning a planar parallelogram in \#86 and a 3D tetrahedron in \#739 and \#756) whose centroid $\mathbf c_1(S)$ is also $\mathbf 0$.}
\label{tab:three-exceptions}
\small
\setlength{\tabcolsep}{4.5pt}
\begin{tabular}{r c c c c c l}
\toprule
Gem & Dudeney & $f(\mathcal G_0)$ & $\vol(\mathcal G_0)$ & $\operatorname{Iso}(\mathcal G_0)$ & $\dim/\vol(\mathcal G_1)$ & First row of standard form $S$\\
\midrule
\#86  & Group VI & $(12,28,18)$ & $209/3$ & $C_2$ ($\overline S=vS$) & $2\mathrm{D}$ ($0$) & $(1,7,14,12)$\\
\#739 & Group VI & $(12,26,16)$ & $287/3$ & $C_2$ ($\overline S=hS$) & $3\mathrm{D}$ ($5/3$) & $(4,14,3,13)$\\
\#756 & Group VI & $(12,24,14)$ & $84$    & $C_2$ ($\overline S=hS$) & $3\mathrm{D}$ ($4/3$) & $(5,2,15,12)$\\
\bottomrule
\end{tabular}
\end{table}

\begin{figure}[htbp]
\centering
\includegraphics[width=\textwidth]{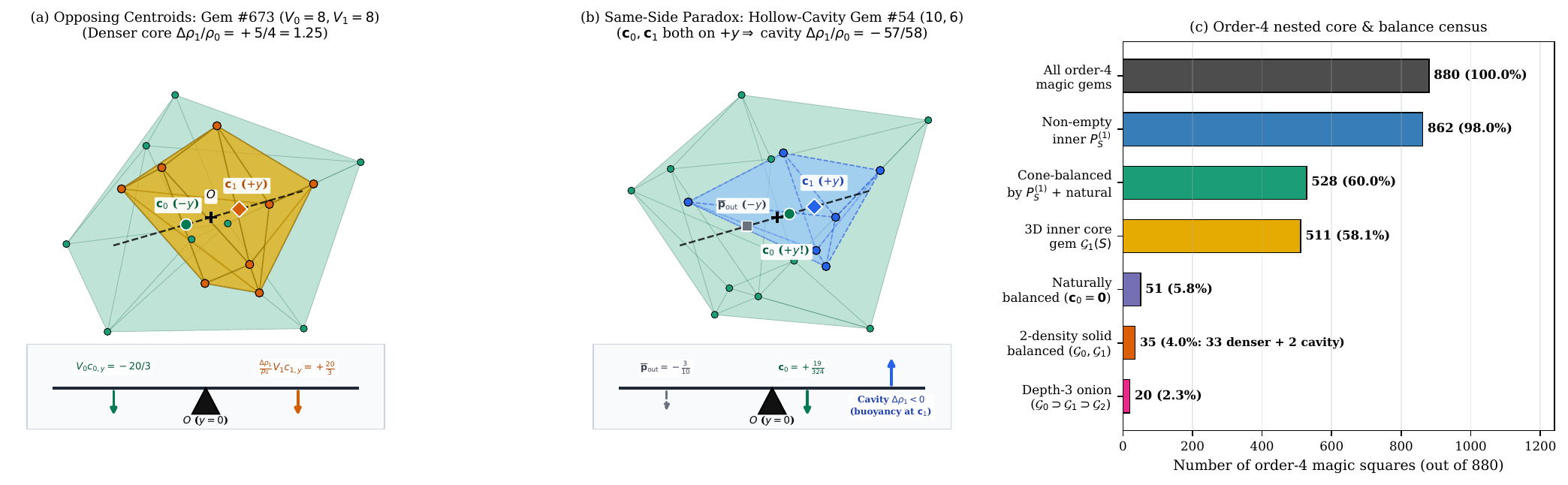}
\caption{Nested ``diamond-in-a-diamond'' convex layers $\mathcal G_0(S)\supset\mathcal G_1(S)$ and centre-of-mass restoration. (a)~Denser-core Gem \#673 ($(V_0,V_1)=(8,8)$) and (b)~hollow-cavity Gem \#54 ($(V_0,V_1)=(10,6)$, $1$-swap sibling of \#45), showing the translucent outer hull $\mathcal G_0(S)$ around the inner core $\mathcal G_1(S)$ with 2D fulcrum lever diagrams below. On \#673, $\mathbf c_0$ (green circle) and $\mathbf c_1$ (gold diamond) lie on opposite sides of $O$, balancing via a heavier inner core ($\Delta\rho_1/\rho_0=+5/4$). On \#54, volumetric simplex weighting pulls $\mathbf c_0=(0,19/324,0)$ onto the $+y$ side alongside $\mathbf c_1=(0,29/68,0)$, so a $98.3\%$ hollow core cavity ($\Delta\rho_1/\rho_0=-57/58$) restores balance at $O$. (c)~Census of nested layer dimensions and solid balance mechanisms across all $880$ order-four magic squares.}
\label{fig:nested-diamonds}
\end{figure}

To counter the centroid drift on the remaining $829$ squares, we peel $P_S$ into nested convex layers \cite{chazelle1985convex}:
\begin{equation}
\label{eq:peeling}
P_S^{(0)} = P_S,\qquad
\mathcal G_k(S) = \conv P_S^{(k)},\qquad
P_S^{(k+1)} = P_S^{(k)} \setminus \vertset\mathcal G_k(S)\quad (k\ge0),
\end{equation}
calling $\mathcal G_0(S)$ the \emph{outer gem} and $\mathcal G_1(S)=\conv P_S^{(1)}\subset\mathcal G_0(S)$ the \emph{inner core gem}—``a diamond inside a diamond'' (\Cref{fig:nested-diamonds}(a,b)).

\begin{proposition}[Vertex-centroid opposition and order-four onion hierarchy]
\label{prop:opposition}
Let $S$ be any arrangement of order $n$ with $V_0=|\vertset\mathcal G_0(S)|<N$ and $V_1=N-V_0$. The centroid $\overline{\mathbf p}_{\mathrm{outer}}=\frac{1}{V_0}\sum_{\mathbf p\in\vertset\mathcal G_0}\mathbf p$ of the outer hull vertices and the centroid $\overline{\mathbf p}_{\mathrm{inner}}=\frac{1}{V_1}\sum_{\mathbf p\in P_S^{(1)}}\mathbf p$ of the inner core points satisfy
\begin{equation}
\label{eq:opposition}
V_0\,\overline{\mathbf p}_{\mathrm{outer}} + V_1\,\overline{\mathbf p}_{\mathrm{inner}} = \mathbf 0,
\end{equation}
pointing in strictly opposite directions across $O$ whenever $\overline{\mathbf p}_{\mathrm{outer}}\ne\mathbf 0$. Across the $880$ magic squares of order~$4$ (\Cref{fig:nested-diamonds}(c)), $18$ squares have $V_0=16$ ($P_S^{(1)}=\emptyset$), $842$ squares have onion depth~$2$ ($14$ singletons, $162$ segments, $175$ planar polygons, and $491$ 3D polyhedra $\mathcal G_1(S)$), and $20$ squares have onion depth~$3$ ($\mathcal G_0\supset\mathcal G_1\supset\mathcal G_2$, comprising $16$ with dimension sequence $(3,3,0)$ and $4$ with $(3,3,1)$), so in total $511$ order-four magic squares possess a three-dimensional inner core gem $\mathcal G_1(S)$.
\end{proposition}

\begin{proof}
Equation~\eqref{eq:opposition} follows immediately from $\sum_{\mathbf p\in P_S}\mathbf p=\mathbf 0$ and $P_S=\vertset\mathcal G_0(S)\sqcup P_S^{(1)}$, and exact convex-hull peeling across all $880$ squares yields the stated layer census.
\end{proof}

\subsection{Restoring Solid Balance via ``Diamond-in-a-Diamond'' Core Inclusions}

Guided by the opposition identity \eqref{eq:opposition}, our first restoration method plays the inner core $\mathcal G_1(S)$ directly against the outer solid hull $\mathcal G_0(S)$. For a two-density composite solid to balance via a single scalar density increment $\Delta\rho_1$, the 3D vectors $\mathbf c_0(S)$ and $\mathbf c_1(S)$ must be collinear—which occurs in Dudeney's Group~VI, where $C_2$ half-turn symmetry locks both centroids onto the same coordinate axis.

\begin{theorem}[Composite ``diamond-in-a-diamond'' solid balancing]
\label{thm:diamond-balance}
Let $S$ be an order-four magic square with $\dim\mathcal G_1(S)=3$ ($511$ squares), and fabricate $\mathcal G_0(S)$ with baseline density $\rho_0>0$ and inner core inclusion $\mathcal G_1(S)\subset\mathcal G_0(S)$ of excess density $\Delta\rho_1=\rho_1-\rho_0$. The composite centre of mass
\[
\mathbf c_{\mathrm{comp}} = \frac{\rho_0\vol(\mathcal G_0)\,\mathbf c_0 + \Delta\rho_1\vol(\mathcal G_1)\,\mathbf c_1}{\rho_0\vol(\mathcal G_0) + \Delta\rho_1\vol(\mathcal G_1)}
\]
vanishes for some $\Delta\rho_1\ne0$ if and only if $\mathbf c_0$ and $\mathbf c_1$ are non-zero and collinear ($\mathbf c_0\times\mathbf c_1=\mathbf 0$), with unique balancing density increment
\begin{equation}
\label{eq:density-ratio}
\frac{\Delta\rho_1}{\rho_0} = -\frac{\vol(\mathcal G_0)\,(\mathbf c_0\cdot\mathbf c_1)}{\vol(\mathcal G_1)\,\|\mathbf c_1\|_2^2}.
\end{equation}
Across the $511$ order-four gems with $\dim\mathcal G_1(S)=3$, $207$ gems have $\mathbf c_1=\mathbf 0$ ($13$ associative gems plus $194$ non-associative gems whose four inner points form a zero-sum tetrahedron) and $254$ gems outside Group~VI have non-collinear solid centroids $\mathbf c_0\times\mathbf c_1\ne\mathbf 0$. Among the remaining $50$ gems (all in Dudeney's Group~VI), $33$ have strictly anti-parallel solid centroids ($\mathbf c_0\cdot\mathbf c_1<0$), balanced by a denser inner core with rational excess density $\Delta\rho_1/\rho_0\in[5/4,\,191/6]$ (\Cref{tab:two-density} and \Cref{tab:all33}), led by Gem \#673 at $\Delta\rho_1/\rho_0=5/4$ (\Cref{fig:hero-squares}(e), \Cref{fig:nested-diamonds}(a)); the other $17$ have parallel solid centroids pointing in the same direction ($\mathbf c_0\cdot\mathbf c_1>0$), because volumetric simplex weighting pulls $\mathbf c_0$ onto the same side of $O$ as $\mathbf c_1$ despite the opposite vertex centroids of \eqref{eq:opposition}, and two of these $17$ gems balance with non-negative core density $\rho_1\in[0,\rho_0)$: square \#123 ($\Delta\rho_1/\rho_0=-3/17$, $\rho_1/\rho_0=14/17$, a lighter core) and square \#54 ($\Delta\rho_1/\rho_0=-57/58$, $\rho_1/\rho_0=1/58$, a nearly hollow core cavity, \Cref{fig:nested-diamonds}(b)), which differs from denser-core square \#45 ($\Delta\rho_1/\rho_0=+85/61$) by the single border swap $(8,9)\leftrightarrow(15,2)$ (\Cref{fig:hero-squares}(c)).

More generally, placing non-negative point masses $m_q\ge0$ at the inner points $q\in P_S^{(1)}$ balances the uniform solid $\mathcal G_0(S)$ at $O$ if and only if $-\mathbf c_0(S)\in\cone(P_S^{(1)})$, which holds for $477$ of the $829$ unbalanced gems ($528/880=60.0\%$ including the $51$ naturally balanced solids), while partitioning $\mathcal G_0(S)$ into $N$ cells of equal mass $M/N$ with centroids at $\mathbf p_{ij}$ (or embedding $N$ identical dense beads at $\mathbf p_{ij}$ in a lightweight shell) yields $\mathbf c_{\mathrm{comp}}=\mathbf 0$ and restores the axisymmetric inertia tensor of \Cref{thm:rotor} for every semi-magic square of any order $n$.
\end{theorem}

\begin{proof}
Setting $\rho_0\vol(\mathcal G_0)\mathbf c_0+\Delta\rho_1\vol(\mathcal G_1)\mathbf c_1=\mathbf 0$ and evaluating the exact rational simplicial moments (\Cref{sec:appendix}) across all $511$ gems with $\dim\mathcal G_1(S)=3$ yields \eqref{eq:density-ratio} and the four-class partition. Existence of non-negative inner point masses $m_q\ge0$ is equivalent to $-\mathbf c_0\in\cone(P_S^{(1)})$ and certified by linear programming over the $829$ unbalanced gems, while the bead-inlaid identity follows directly from $\sum_{i,j}\mathbf p_{ij}=\mathbf 0$ and \Cref{thm:rotor}.
\end{proof}

\begin{table}[htbp]
\centering
\caption{The eight order-four Magic Gems with the smallest balancing excess density ratios $\Delta\rho_1/\rho_0\in[5/4,40/19]\approx[1.25,2.11]$ in a two-density ``diamond-in-a-diamond'' solid inclusion $(\mathcal G_0(S),\mathcal G_1(S))$ (\Cref{thm:diamond-balance}), all achievable with standard dual-material 3D printing. All $33$ denser-core gems belong to Dudeney's Group~VI; the full list of $33$ is given in \Cref{tab:all33}.}
\label{tab:two-density}
\small
\setlength{\tabcolsep}{5.0pt}
\begin{tabular}{r c c c c c c c}
\toprule
Gem & $(V_0,V_1)$ & $\vol(\mathcal G_0)$ & $\vol(\mathcal G_1)$ & Axial $c_0$ & Axial $c_1$ & Exact $\Delta\rho_1/\rho_0$ & Decimal $\rho_1/\rho_0$\\
\midrule
\#673 & $(8,8)$  & $257/3$ & $56/3$ & $-20/257$  & $+2/7$    & $5/4$   & $2.250$\\
\#502 & $(10,6)$ & $271/3$ & $6$    & $+2/271$   & $-1/12$   & $4/3$   & $2.333$\\
\#45  & $(10,6)$ & $214/3$ & $6$    & $-85/856$  & $+61/72$  & $85/61$ & $2.393$\\
\#50  & $(10,6)$ & $254/3$ & $22/3$ & $-7/127$   & $+19/44$  & $28/19$ & $2.474$\\
\#307 & $(10,6)$ & $277/3$ & $22/3$ & $+23/1108$ & $-15/88$  & $23/15$ & $2.533$\\
\#49  & $(12,4)$ & $81$    & $2$    & $-2/81$    & $+1/2$    & $2$     & $3.000$\\
\#38  & $(10,6)$ & $263/3$ & $22/3$ & $+31/1052$ & $-15/88$  & $31/15$ & $3.067$\\
\#47  & $(10,6)$ & $247/3$ & $22/3$ & $-20/247$  & $+19/44$  & $40/19$ & $3.105$\\
\bottomrule
\end{tabular}
\end{table}

\subsection{Overcoming the Density Barrier: Complete \texorpdfstring{$K_N$}{K\_N} Networks, Nested Wireframes, and Papercraft Shells}
\label{subsec:wireframes-nets}

While the eight gems in \Cref{tab:two-density} have modest excess density ratios $\Delta\rho_1/\rho_0\in[5/4,40/19]$ that print directly on a dual-extruder 3D printer (such as standard PLA around a tungsten- or bronze-filled PLA core), the remaining $25$ denser-core gems in \Cref{tab:all33} have small core volumes $\vol(\mathcal G_1)$ requiring $\Delta\rho_1/\rho_0\in[47/12,191/6]$, and $194$ gems have a four-vertex tetrahedral core whose solid centroid $\mathbf c_{\mathrm{solid}}(\mathcal G_1)=\frac14\sum\mathbf p_q=\mathbf 0$ provides no lever arm at all.

\begin{figure}[htbp]
\centering
\includegraphics[width=\textwidth]{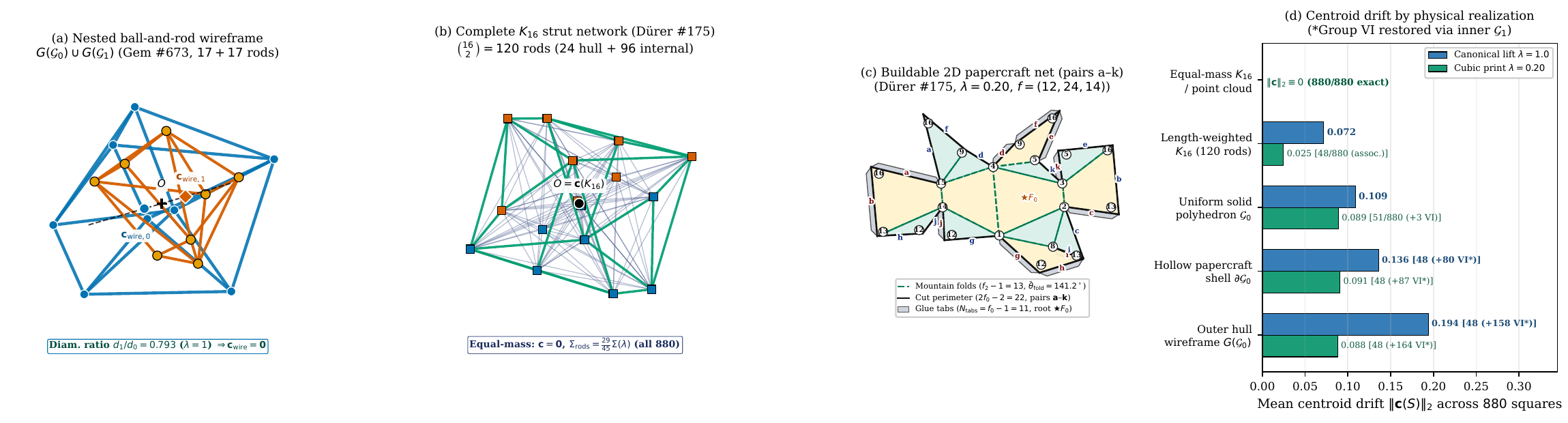}
\caption{Wireframes, complete $K_{16}$ strut networks, and unfolded papercraft cutout nets of Magic Gems \cite{mathewson2025magicgemweb}. (a)~Nested ball-and-rod wireframe $G(\mathcal G_0)\cup G(\mathcal G_1)$ of Gem \#673 ($17$ outer blue rods and $17$ inner orange rods), balanced at rod-diameter ratio $d_1/d_0=0.793$ ($\lambda=1$). (b)~Complete network $K_{16}$ of D\"urer's Gem \#175 with all $\binom{16}{2}=120$ rods ($24$ hull edges in green, $96$ internal cross-struts in slate), which satisfies $\mathbf c(K_{16})=\mathbf 0$ and $\Sigma_{\mathrm{rods}}=\frac{29}{45}\Sigma(\lambda)$ for equal-mass rods across all $880$ squares. (c)~Exact 2D unfolded papercraft net of D\"urer's Gem \#175 ($\lambda=0.20$, $f=(12,24,14)$) rooted at $\bigstar F_0$, with $f_2-1=13$ dashed fold creases and $N_{\mathrm{tabs}}=f_0-1=11$ glue tabs paired by labels \texttt{a}--\texttt{k}. (d)~Mean centroid drift $\|\mathbf c(S)\|_2$ across all $880$ squares for the five physical realizations at $\lambda=1.0$ and $\lambda=0.20$.}
\label{fig:wireframes-nets}
\end{figure}

Replacing solid volume with the skeletal and thin-shell architectures of \texttt{magicGemWeb} \cite{mathewson2025magicgemweb} (\Cref{fig:wireframes-nets}) overcomes both barriers. Connecting all $\binom{N}{2}$ pairs of points as a complete strut network $K_N\bigl(P_S(\lambda)\bigr)$ (\texttt{magicGem\_v2.stl}, \Cref{fig:wireframes-nets}(b)) restores equal participation of all $N$ vertices; building a nested ball-and-rod wireframe $G(\mathcal G_0)\cup G(\mathcal G_1)$ (\texttt{magicGem\_v1.stl}, \Cref{fig:wireframes-nets}(a)) replaces volumetric scaling with quadratic rod-diameter scaling ($\rho_{\mathrm{lin}}\propto d^2$) and gives four-vertex tetrahedral cores a non-zero edge-length centroid $\mathbf c_{\mathrm{wire}}(\mathcal G_1)\ne\mathbf 0$; and unfolding $\partial\mathcal G_\lambda(S)$ along a dual spanning tree produces a flat printable papercraft net (\Cref{fig:wireframes-nets}(c)).

\begin{theorem}[Complete $K_N$ strut networks and two-gauge nested wireframes]
\label{thm:complete-network}
Let $S$ be an arrangement of $\{1,\dots,N\}$, $N=n^2$ ($n\ge2$), with $\sum_{u=1}^N\mathbf p_u(\lambda)=\mathbf 0$ and empirical covariance matrix $\Sigma(\lambda)=\frac1N\sum_{u=1}^N\mathbf p_u(\lambda)\mathbf p_u(\lambda)^{\top}$ in the slender-rod limit $d/L_{uv}\to 0^+$. If each of the $\binom{N}{2}$ rods connecting $\mathbf p_u(\lambda)$ and $\mathbf p_v(\lambda)$ ($1\le u<v\le N$) carries equal mass $m_{\mathrm{edge}}=M_{\mathrm{rods}}/\binom{N}{2}$ (with optional equal node masses $m_{\mathrm{node}}=M_{\mathrm{nodes}}/N$), then $(N-1)$-regularity of $K_N$ places the centre of mass at $O$ for every arrangement $S$:
\begin{equation}
\label{eq:kN-com}
\mathbf c\bigl(K_N(\lambda)\bigr) = \frac{1}{\binom{N}{2}}\sum_{1\le u<v\le N}\frac{\mathbf p_u(\lambda)+\mathbf p_v(\lambda)}{2} = \frac{1}{N}\sum_{u=1}^N\mathbf p_u(\lambda) = \mathbf 0,
\end{equation}
and the continuous second-moment covariance matrix of the $\binom{N}{2}$ rods is a universal scalar multiple of $\Sigma(\lambda)$:
\begin{equation}
\label{eq:kN-cov}
\Sigma_{\mathrm{rods}}^{(=\text{mass})}\bigl(K_N(\lambda)\bigr)
= \frac{2N-3}{3(N-1)}\,\Sigma(\lambda)
\qquad\Bigl(=\frac{5}{8}\,\Sigma(\lambda)\text{ at }n=3,\quad \frac{29}{45}\,\Sigma(\lambda)\text{ at }n=4\Bigr),
\end{equation}
so that for every semi-magic square $S$, the equal-mass complete network $K_N\bigl(P_S(\lambda)\bigr)$ shares the exact axisymmetric inertia ratio $I_\perp(\lambda)/I_{zz}=\frac12\bigl(1+\lambda^2(n^2+1)\bigr)$ and spherical-top scale $\lambda_{\mathrm{iso}}=(n^2+1)^{-1/2}$ of \Cref{thm:rotor}.

When a subgraph $G=(V_G,E_G)$ is instead printed with uniform rod diameter $d$ so that each rod has mass proportional to its length $L_{uv}(\lambda)=\|\mathbf p_u(\lambda)-\mathbf p_v(\lambda)\|_2$, both $K_{16}$ and the outer wireframe $G(\mathcal G_0)$ have $\mathbf c_{\mathrm{len}}=\mathbf 0$ if and only if $S$ is one of the $48$ associative squares, while passing from $G(\mathcal G_0)$ to the $120$-rod complete network $K_{16}$ reduces mean centroid drift $\|\mathbf c_{\mathrm{len}}\|_2$ across $\mathcal M_4$ by $2.7\times$ at $\lambda=1$ (from $0.194$ to $0.072$) and by $3.6\times$ at $\lambda=0.20$ (from $0.088$ to $0.025$; \Cref{fig:wireframes-nets}(d)). Furthermore, for the $848$ order-four gems with $\dim\mathcal G_1(S)\ge1$, printing the outer edges $E(\mathcal G_0)$ with rod diameter $d_0$ and the disjoint inner core edges $E(\mathcal G_1)$ with rod diameter $d_1$ restores $\mathbf c_{\mathrm{wire}}(\mathcal G_0\cup\mathcal G_1)=\mathbf 0$ whenever $\mathbf c_{\mathrm{wire}}(\mathcal G_0)$ and $\mathbf c_{\mathrm{wire}}(\mathcal G_1)$ are strictly anti-parallel, at the unique rod diameter ratio
\begin{equation}
\label{eq:rod-diam-ratio}
\frac{d_1}{d_0} = \sqrt{\frac{\rho_{\mathrm{lin},1}}{\rho_{\mathrm{lin},0}}}
= \sqrt{-\frac{L\bigl(G(\mathcal G_0)\bigr)\,\bigl(\mathbf c_{\mathrm{wire}}(\mathcal G_0)\cdot\mathbf c_{\mathrm{wire}}(\mathcal G_1)\bigr)}{L\bigl(G(\mathcal G_1)\bigr)\,\|\mathbf c_{\mathrm{wire}}(\mathcal G_1)\|_2^2}},
\end{equation}
balancing $158$ Group~VI gems at $\lambda=1$ (such as Gem \#673 at $d_1/d_0=0.793$ and Gem \#45 at $d_1/d_0=0.795$) and $164$ Group~VI gems at $\lambda=0.20$.
\end{theorem}

\begin{proof}
Integrating $\mathbf x(t)=(1-t)\mathbf p_u+t\mathbf p_v$ over $t\in[0,1]$ gives midpoint $\frac12(\mathbf p_u+\mathbf p_v)$ and second moment $\frac13(\mathbf p_u\mathbf p_u^{\top}+\mathbf p_v\mathbf p_v^{\top})+\frac16(\mathbf p_u\mathbf p_v^{\top}+\mathbf p_v\mathbf p_u^{\top})$. Since each vertex belongs to $N-1$ edges of $K_N$ and $\sum_{u<v}(\mathbf p_u\mathbf p_v^{\top}+\mathbf p_v\mathbf p_u^{\top})=-N\Sigma(\lambda)$, averaging over $\binom{N}{2}$ edges yields \eqref{eq:kN-com}--\eqref{eq:kN-cov}. Exact evaluation of the edge lengths and \eqref{eq:rod-diam-ratio} across all $880$ squares establishes the length-weighted and nested wireframe censuses (\Cref{sec:appendix}).
\end{proof}

\begin{proposition}[Papercraft cutout nets, glue-tab invariant, and nested thin-shell balance]
\label{prop:papercraft-net}
Let $\mathcal G_\lambda(S)$ be any three-dimensional Magic Gem with $f$-vector $(f_0,f_1,f_2)$. Unfolding $\partial\mathcal G_\lambda(S)$ into a connected planar net along any spanning tree $T^*$ of the dual facet graph $G^*(\mathcal G)$ \cite{durer1525underweysung,shephard1975convex,demaine2007geometric,ghomi2014durer} partitions the $f_1$ edges into $f_2-1$ interior fold creases and $N_{\mathrm{tabs}}=f_1-(f_2-1)$ cut edges joined by trapezoidal glue tabs (paired by labels \texttt{a}--\texttt{k} in \Cref{fig:wireframes-nets}(c)). Because $\partial\mathcal G_\lambda(S)\cong S^2$, the cut edges $E\setminus E(T^*)$ form a spanning tree $T$ of the primal $1$-skeleton $G(\mathcal G_\lambda)$ on the $f_0$ hull vertices, so the number of glue tabs depends solely on the hull vertex count $f_0$:
\begin{equation}
\label{eq:glue-tabs}
N_{\mathrm{tabs}} = f_1 - (f_2 - 1) = |E(T)| = f_0 - 1,
\end{equation}
taking every integer value in $\{7,\dots,15\}$ across $\mathcal M_4$ ($11$ tabs for D\"urer's \#175). Under the deterministic breadth-first dual spanning tree $T^*$ rooted at the largest-area facet $F_0$ ($\bigstar F_0$ in \Cref{fig:wireframes-nets}(c)), all $880$ order-four Magic Gems ($100\%$) unfold into strictly non-overlapping 2D planar nets at both $\lambda=1$ and $\lambda=0.20$. Moreover, a uniform boundary shell $\partial\mathcal G_\lambda(S)$ of areal density $\sigma_0>0$ has $\mathbf c_{\mathrm{shell}}(\mathcal G_0)=\mathbf 0$ on the $48$ associative squares, and enclosing a nested inner shell $\partial\mathcal G_1(S)$ ($\dim\mathcal G_1=3$) of areal density $\sigma_1$ restores $\mathbf c_{\mathrm{shell}}(\partial\mathcal G_0\cup\partial\mathcal G_1)=\mathbf 0$ for $80$ Group~VI gems at $\lambda=1$ ($87$ at $\lambda=0.20$).
\end{proposition}

\begin{proof}
By planar graph duality on $\partial\mathcal G_\lambda(S)\cong S^2$, the complement of a dual spanning tree $T^*\subset G^*(\mathcal G_\lambda)$ ($|E(T^*)|=f_2-1$) is a primal spanning tree $T\subset G(\mathcal G_\lambda)$ on the $f_0$ vertices, giving $N_{\mathrm{tabs}}=|E(T)|=f_0-1$ (invariant under facet triangulation). Exact 2D polygon intersection testing and triangular surface moments certify non-overlap and the shell censuses across all $880$ gems (\Cref{sec:appendix}).
\end{proof}

\begin{remark}[Physical realization of finite-diameter struts and core suspension]
\label{rem:finite-diameter}
In a physical 3D print of finite rod diameter $d>0$, each cylindrical strut of length $L_{uv}$ and direction $\hat{\mathbf e}_{uv}$ retains midpoint centroid $\frac12(\mathbf p_u+\mathbf p_v)$ and adds the radial inertia correction $m_{uv}\frac{d^2}{16}(I_3+\hat{\mathbf e}_{uv}\hat{\mathbf e}_{uv}^{\top})$ up to $O(d^3)$ hub-overlap volumes. Suspending the inner core wireframe $G(\mathcal G_1)$ inside $G(\mathcal G_0)$ on a Group~VI gem via $C_2$-symmetric pairs of equal-length struts places the suspension centroid at $O$, leaving the balancing diameter ratio $d_1/d_0$ in \eqref{eq:rod-diam-ratio} invariant.
\end{remark}

Having restored internal balance to individual polyhedra, we now remove the central pivot $O$ in \Cref{sec:stacking} and place Magic Gems directly on a flat horizontal table and into multi-block vertical stacks under gravity.

\section{Freestanding Gravity and Mirror-Twin Towers}
\label{sec:stacking}

Our fourth and final test removes the central pivot and places Magic Gems directly onto a flat horizontal table under gravity—first as individual resting blocks, and then stacked face-to-face into multi-block vertical columns. Because the hull vertices sit at unequal heights, the boundary facets of $\mathcal G_\lambda(S)$ are slanted polygons whose outward normals rarely oppose each other across the gem. We therefore test when a slanted facet $F\subset\partial\mathcal G_\lambda(S)$ supports the gem's centre of mass without tipping, and how pairing a gem with an inverted partner restores a straight vertical column when identical copies curve off-axis.

\subsection{Single-Block Table Equilibrium: Cubic Aspect Scaling Stabilizes \texorpdfstring{$74.5\%$}{74.5\%} of Facets}

Let $F\subset\partial\mathcal G_\lambda(S)$ be a facet with outward unit normal $\hat{\mathbf n}_F(\lambda)$. When $\mathcal G_\lambda(S)$ rests with $F$ on a horizontal table (so $\hat{\mathbf n}_F(\lambda)$ points straight down), gravity acts vertically through the centre of mass $\mathbf c(\lambda)$. Let $\pi_F(\mathbf c(\lambda))$ denote the orthogonal projection of $\mathbf c(\lambda)$ onto the affine plane of $F$, and let $h_F(\lambda)=(\mathbf p_F-\mathbf c(\lambda))\cdot(-\hat{\mathbf n}_F(\lambda))>0$ be the height of the centre of mass above the table.

\begin{definition}[Static stability margin and critical tipping angle]
\label{def:stability}
For a facet $F\subset\partial\mathcal G_\lambda(S)$, the \emph{stability margin} $d_{\mathrm{stab}}(S,F;\lambda)$ is the signed Euclidean distance in the plane of $F$ from $\pi_F(\mathbf c(\lambda))$ to the boundary $\partial F$ of the facet polygon (positive when $\pi_F(\mathbf c(\lambda))$ lies in the strict relative interior of $F$, and negative outside). Facet $F$ is \emph{statically stable} if and only if $d_{\mathrm{stab}}(S,F;\lambda)>0$, with \emph{critical tipping angle}
\begin{equation}
\label{eq:tip-angle}
\theta_{\mathrm{tip}}(S,F;\lambda) = \arctan\!\left(\frac{d_{\mathrm{stab}}(S,F;\lambda)}{h_F(\lambda)}\right).
\end{equation}
\end{definition}

\begin{proposition}[Order-four resting stability census]
\label{prop:stability-census}
Across all $880$ magic squares of order~$4$, the gems $\mathcal G_\lambda(S)$ possess a total of $14{,}824$ facets (mean $16.85$ facets per gem). At the tall canonical lift $\lambda=1$ (where the vertical height span $15$ is five times the $3\times3$ floor span), the high centre of mass $\mathbf c_0(S)$ projects outside $57.0\%$ of all facet polygons: $6{,}378$ of the $14{,}824$ facets ($43.0\%$) are statically stable, with every gem possessing between $4$ and $12$ stable resting facets (mean $7.25$; using the balanced point centroid $O$ gives $6{,}340$ stable facets, $42.8\%$). Compressing the vertical lift to the cubic 3D-print scale $\lambda_{\mathrm{cube}}=1/5=0.20$ brings the centre of mass closer to the slanted upper and lower facets so that $11{,}044$ of the $14{,}824$ facets ($74.5\%$) become statically stable, with a mean critical tipping angle of $17.8^\circ$ and a maximum of $47.3^\circ$ (\Cref{fig:mirror-stacking}(d)).
\end{proposition}

\begin{proof}
For each facet $F\subset\partial\mathcal G_\lambda(S)$ with exact integer plane equation $aX+bY+cZ+d=0$ in doubled coordinates, rescaling $z\mapsto\lambda z$ transforms the outward normal to $\mathbf n_F(\lambda)\propto(a\lambda,b\lambda,c)$ and the solid centroid to $\mathbf c_0(\lambda)=(c_{0,x},c_{0,y},\lambda c_{0,z})$. Projecting $\mathbf c_0(\lambda)$ onto the plane of $F$ and evaluating $d_{\mathrm{stab}}(S,F;\lambda)$ and \eqref{eq:tip-angle} across all $14{,}824$ facets yields the stated counts and angles.
\end{proof}

\subsection{Why Same-Gem Stacks Wedge Apart, and How Mirror-Twin Inversion Restores Verticality}

\begin{figure}[htbp]
\centering
\includegraphics[width=\textwidth]{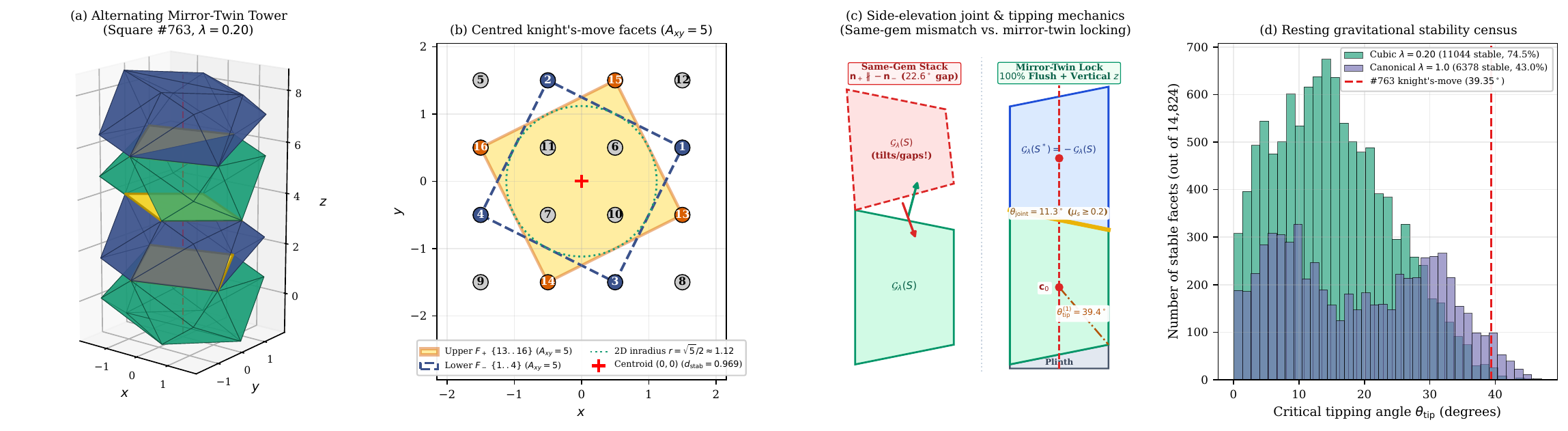}
\caption{Mirror-twin facet locking and periodic vertical towers. (a)~Four-block alternating tower of square \#763 at $\lambda=0.20$, alternating $\mathcal G_\lambda(S)$ (green) with its mirror twin $\mathcal G_\lambda(r^2\overline S)=-\mathcal G_\lambda(S)$ (blue) along the golden knight's-move contact facets $\{13,14,15,16\}$ and $\{1,2,3,4\}$ to keep the column's centre of mass on the $z$-axis (red dashed). (b)~Top-down $(x,y)$ footprint of square \#763 showing the dual knight's-move squares of area $A_{xy}=5$ and 2D inradius $r=\sqrt{5}/2\approx1.12$ ($d_{\mathrm{stab}}=0.969$ in the tilted plane at $\lambda=0.20$). (c)~Side-elevation schematic comparing the $2\theta_{\mathrm{joint}}=22.6^\circ$ wedge mismatch of same-gem stacking against flush mirror-twin locking ($\mathbf n_+^*=-\mathbf n_+$) on an $11.3^\circ$ base plinth. (d)~Census of critical tipping angles $\theta_{\mathrm{tip}}$ over all $14{,}824$ facets at $\lambda=0.20$ and $\lambda=1.0$.}
\label{fig:mirror-stacking}
\end{figure}

Suppose we now stack multiple blocks vertically along an upper facet $F_+\subset\mathcal G_\lambda(S)$ (outward normal $+\mathbf n_+$, $n_{+,z}>0$) and a lower facet $F_-\subset\mathcal G_\lambda(S)$ (outward normal $+\mathbf n_-$, $n_{-,z}<0$). Stacking translated copies of the \emph{same} gem $\mathcal G_\lambda(S)$ without rotation requires parallel opposite normals $\mathbf n_-=-\mathbf n_+$, which fails on $756$ of the $880$ gems ($85.9\%$); when $\mathbf n_-\ne-\mathbf n_+$, identical copies meet at an open wedge angle that bends the column away from the vertical (\Cref{fig:mirror-stacking}(c), left).

Equilibrium is restored by pairing a gem with its \emph{rotational complement (mirror twin)}
\begin{equation}
\label{eq:mirror-twin}
S^* := r^2\overline S = (N+1)\mathbf 1 - r^2 S,
\end{equation}
which inverts the polyhedron through the origin in $\R^3$ and pairs every slanted upper facet with an equal and opposite lower facet (\Cref{fig:mirror-stacking}(a,c)).

\begin{theorem}[Mirror-Twin Inversion and Zero-Rotation Alternating Towers]
\label{thm:mirror-tower}
Let $S$ be any arrangement of order $n$, let $S^*=r^2\overline S$ be its half-turned complement, and let $\lambda>0$. In $\R^3$, the lifted point clouds and convex hulls satisfy the central inversion identity
\begin{equation}
\label{eq:central-inv}
P_{S^*}(\lambda) = -P_S(\lambda),\qquad
\mathcal G_\lambda(S^*) = -\mathcal G_\lambda(S),
\end{equation}
so the map $F\mapsto F^*:=-F$ is an area-preserving, normal-reversing bijection between the upper facets $F_+\subset\mathcal G_\lambda(S)$ (outward normal $+\mathbf n_+$) and the lower facets $F_+^*\subset\mathcal G_\lambda(S^*)$ (outward normal $-\mathbf n_+$), and between the lower facets $F_-\subset\mathcal G_\lambda(S)$ (outward normal $+\mathbf n_-$) and the upper facets $F_-^*=-F_-\subset\mathcal G_\lambda(S^*)$ (outward normal $-\mathbf n_-$).

Consequently, given any upper facet $F_+\subset\mathcal G_\lambda(S)$ and any lower facet $F_-\subset\mathcal G_\lambda(S)$ (not necessarily parallel) that are centrally symmetric polygons in $\R^3$ with centroids $\mathbf c(F_+)$ and $\mathbf c(F_-)$, the alternating sequence of unrotated blocks $B_0,B_1,B_2,\dots$ defined by
\begin{equation}
\label{eq:tower-blocks}
B_{2m} = \mathcal G_\lambda(S) + m\mathbf T,\qquad
B_{2m+1} = -\mathcal G_\lambda(S) + 2\mathbf c(F_+) + m\mathbf T,\qquad
\mathbf T = 2\bigl(\mathbf c(F_+) - \mathbf c(F_-)\bigr),
\end{equation}
meets with $100\%$ polygon and vertex overlap at every odd joint $B_{2m}\cap B_{2m+1}$ (where the upper facet $F_+ + m\mathbf T$ of $B_{2m}$ coincides identically with the lower facet $-F_+ + 2\mathbf c(F_+) + m\mathbf T$ of $B_{2m+1}$ with opposite outward normals $\pm\mathbf n_+$) and at every even joint $B_{2m-1}\cap B_{2m}$ (where the upper facet $-F_- + 2\mathbf c(F_+) + (m-1)\mathbf T$ of $B_{2m-1}$ coincides identically with the lower facet $F_- + m\mathbf T$ of $B_{2m}$ with opposite outward normals $\mp\mathbf n_-$). Every two-block repeating unit $B_{2m}\cup B_{2m+1}$ is centrally symmetric about $\mathbf c(F_+)+m\mathbf T$, and whenever $\mathbf c(F_+)=(0,0,z_+)$ and $\mathbf c(F_-)=(0,0,z_-)$ both lie on the vertical $z$-axis ($z_+>0>z_-$), the period vector $\mathbf T=(0,0,2(z_+-z_-))$ is purely vertical and places the centre of mass of any $2M$-block tower on the $z$-axis regardless of whether $\mathbf c_0(S)=\mathbf 0$ (\Cref{fig:mirror-stacking}(a,c)).
\end{theorem}

\begin{proof}
In centred coordinates \eqref{eq:coords}, $r^2$ sends $(x,y)\mapsto(-x,-y)$ and complementation sends $z\mapsto-z$, so $r^2\overline S$ lifts to $-\mathbf p_{ij}(\lambda)$, proving \eqref{eq:central-inv}. Central symmetry of $F_\pm$ about $\mathbf c(F_\pm)$ gives $-F_\pm + 2\mathbf c(F_\pm)=F_\pm$; translating by $m\mathbf T$ establishes identical facet coincidence at all odd and even joints, while central inversion about $\mathbf c(F_+)+m\mathbf T$ exchanges $B_{2m}$ and $B_{2m+1}$ and locks the two-block centroid onto the $z$-axis.
\end{proof}

\subsection{Chiral Knight's-Move Towers: Exemplar Square \#763}

\begin{proposition}[Vertically centered tower census and the knight's-move square \#763]
\label{prop:knights-move}
Across the $880$ magic squares of order~$4$, exactly $48$ squares possess both an upper facet $F_+$ and a lower facet $F_-$ that are centrally symmetric polygons whose centroids lie on the vertical axis $(0,0,z_\pm)$. The maximum projected horizontal contact area among these $48$ squares is $A_{xy}=5$ (spanning $5/9\approx55.6\%$ of the $3\times3$ grid footprint), achieved by $17$ squares (\Cref{tab:knights-move}; e.g.\ squares \#603, \#763, and \#774) in which $\pi_{xy}(F_+)$ is the centred $2\times1$ knight's-move square with vertices $(\pm\frac32,\mp\frac12)$ and $(\pm\frac12,\pm\tfrac32)$ (\Cref{fig:mirror-stacking}(b)). Only two of the $17$ (\#774 and \#792) have parallel top and bottom facets ($\mathbf n_-=-\mathbf n_+$); the remaining $15$ require mirror-twin alternation.

In particular, square \#763 (Dudeney Group~VI, $\overline S=hS$, so $S^*=r^2\overline S=vS$ is its vertical reflection) is obtained from naturally balanced solid \#756 ($\mathbf c_0=\mathbf 0$, \Cref{tab:three-exceptions}) by swapping the single complementary pair $(4,13)\leftrightarrow(16,1)$ in rows~$2$--$3$ (\Cref{fig:hero-squares}(f)):
\[
S_{\#763} = \begin{pmatrix}
5 & 2 & \mathbf{15} & 12\\
\mathbf{16} & 11 & 6 & 1\\
4 & 7 & 10 & \mathbf{13}\\
9 & \mathbf{14} & 3 & 8
\end{pmatrix},\qquad
S^*_{\#763} = vS_{\#763} = \begin{pmatrix}
9 & 14 & 3 & 8\\
4 & 7 & 10 & 13\\
16 & 11 & 6 & 1\\
5 & 2 & 15 & 12
\end{pmatrix}.
\]
This swap lifts $\{13,14,15,16\}$ onto the Pythagorean upper knight's-move plane $4x-3y+5z=30$ (centroid $(0,0,6\lambda)$) and $\{1,2,3,4\}$ onto the lower plane $-4x-3y-5z=30$ (centroid $(0,0,-6\lambda)$), while shifting the solid centroid slightly to $\mathbf c_0=(0,-5/126,0)$ (which cancels pairwise in $B_{2m}\cup B_{2m+1}$). Both facets have 3D area $5\sqrt{1+\lambda^2}$, stability margin $d_{\mathrm{stab}}=0.969$, and single-block tipping angle $\theta_{\mathrm{tip}}^{(1)}=39.35^\circ$ at $\lambda=0.20$. Because $\mathbf n_\pm(\lambda)\propto(\pm4\lambda,-3\lambda,\pm5)$ makes angle $\theta_{\mathrm{joint}}(\lambda)=\arctan(\lambda)$ with $\hat{\mathbf e}_z$ ($\approx11.31^\circ$ at $\lambda=0.20$), any material with static friction $\mu_s>\lambda=0.20$ holds the unbonded column at level ($\theta_{\mathrm{table}}=0$), while tilting an $M$-block column of centre-of-mass height $H_M$ is governed by $\min\bigl(\arctan(\mu_s)-\theta_{\mathrm{joint}}(\lambda),\,\arctan(d_{\mathrm{stab}}/H_M)\bigr)$ (\Cref{fig:mirror-stacking}(c)).
\end{proposition}

\begin{proof}
Testing each facet $F\subset\partial\mathcal G_1(S)$ across the $880$ gems for central symmetry ($F=2\mathbf c(F)-F$ with $c_x(F)=c_y(F)=0$) identifies the $48$ gems and the $17$ gems of projected shoelace area $A_{xy}=5$. Substituting the lifted coordinates of $\{13,14,15,16\}$ and $\{1,2,3,4\}$ on \#763 into $4x-3y+5(z/\lambda)=30$ and $-4x-3y-5(z/\lambda)=30$ verifies the facet planes, normals $\mathbf n_\pm(\lambda)\propto(\pm4\lambda,-3\lambda,\pm5)$, surface area $5\sqrt{1+\lambda^2}$, and joint tilt $\tan\theta_{\mathrm{joint}}(\lambda)=\lambda$.
\end{proof}

\begin{remark}[Three-dimensional chirality and why a single mold cannot build the \#763 tower]
\label{rem:chirality}
Because square \#763 satisfies $\overline S=hS$ ($z(-x,y)=-z(x,y)$), its 3D gem $\mathcal G_\lambda(S_{\#763})$ is invariant under the $180^\circ$ proper rotation around the $y$-axis, $R_{y,\pi}(x,y,z)=(-x,y,-z)\in SO(3)$, which exchanges $F_+\colon 4x-3y+5z=30$ and $F_-\colon -4x-3y-5z=30$ and returns the exact same block in its original orientation. Flipping a second copy of $\mathcal G_\lambda(S_{\#763})$ upside-down around the $x$-axis via $R_{x,\pi}(x,y,z)=(x,-y,-z)$ maps $F_+$ to a lower facet with outward normal $(+4\lambda,+3\lambda,-5)\ne-\mathbf n_+=(-4\lambda,+3\lambda,-5)$, failing to mate flush with $F_+$. More fundamentally, $\operatorname{Aut}_{O(3)}\bigl(\mathcal G_\lambda(S_{\#763})\bigr)=\{I_3,R_{y,\pi}\}\cong C_2\subset SO(3)$ contains no improper reflections: $\mathcal G_\lambda(S_{\#763})$ and $\mathcal G_\lambda(S^*_{\#763})=-\mathcal G_\lambda(S_{\#763})$ are \emph{chiral 3D enantiomorphs} in disjoint proper rigid-motion orbits, so building the tower of \Cref{fig:mirror-stacking}(a) requires two opposite-handed molds.
\end{remark}

\begin{table}[htbp]
\centering
\caption{Representative order-four Magic Gems from the $17$ classes (\#603, \#763, \#768, \#774, \#779, \#785, \#788, \#791, \#792, \#818, \#824, \#828, \#836, \#837, \#839, \#842, \#844) achieving the maximum projected contact area $A_{xy}=5$ via centred knight's-move square facets $F_+$ and $F_-$ (\Cref{prop:knights-move}). Plane equations $(a,b,c,d)$ are given in doubled coordinates $aX+bY+cZ+d=0$. ``Same-gem $\parallel$'' indicates whether $\mathbf n_-=-\mathbf n_+$; stability margins $d_{\mathrm{stab}}$ and single-block critical tipping angles $\theta_{\mathrm{tip}}^{(1)}$ (\Cref{def:stability}) are reported for the upper facet at $\lambda=0.20$.}
\label{tab:knights-move}
\small
\setlength{\tabcolsep}{5.0pt}
\begin{tabular}{r c c c c c c}
\toprule
Gem & Group & Upper $F_+$ plane & Lower $F_-$ plane & Same-gem $\parallel$ & $d_{\mathrm{stab}}$ ($\lambda=0.2$) & $\theta_{\mathrm{tip}}^{(1)}$ ($\lambda=0.2$)\\
\midrule
\#603 & VI & $(2,-9,5,-50)$ & $(-2,-9,-5,-50)$ & No  & $1.020$ & $45.39^\circ$\\
\#763 & VI & $(4,-3,5,-60)$ & $(-4,-3,-5,-60)$ & No  & $0.969$ & $39.35^\circ$\\
\#768 & II & $(0,-1,1,-12)$ & $(-4,-3,-5,-60)$ & No  & $0.973$ & $39.24^\circ$\\
\#774 & VI & $(1,0,1,-12)$  & $(-1,0,-1,-12)$  & Yes & $0.931$ & $38.35^\circ$\\
\#791 & VI & $(-3,-4,5,-60)$& $(-3,4,-5,-60)$  & No  & $0.969$ & $39.35^\circ$\\
\#792 & VI & $(0,-1,1,-12)$ & $(0,1,-1,-12)$   & Yes & $0.931$ & $38.35^\circ$\\
\#837 & VI & $(4,-3,5,-60)$ & $(4,3,-5,-60)$   & No  & $0.962$ & $38.94^\circ$\\
\bottomrule
\end{tabular}
\end{table}

\section{Synthesis and Open Questions}
\label{sec:discussion}

Stepping back from the individual theorems, the mechanical tests of Sections~\ref{sec:vectors}--\ref{sec:stacking} tell a unified story about how combinatorial \emph{balance in inequality} behaves when lifted from a flat grid of numbers into three-dimensional physical space (\Cref{tab:synthesis}). Across all four physical regimes—radial spokes, spinning rotors, convex polyhedra, and stacked columns—the $880$ order-four magic squares stratify into three distinct mechanical tiers determined by how their symmetries interact with nonlinearity.

At the broadest level, whenever a physical observable is linear or symmetric-quadratic in the coordinates $(x_j,y_i,z_{ij})$, the row and column sum constraints neutralize entry inequality universally across all $880$ squares ($100\%$). This universal tier encompasses natural 3D and 2D spoke balance ($\beta=1,\alpha=1$), cubic vertical spoke vanishing ($(\mathbf V_3)_z=0$), existence of positive unit-sphere counterweights ($\mathbf w^*>\mathbf 0$), axisymmetric precession and the $\lambda_{\mathrm{iso}}=1/\sqrt{17}$ spherical-top transition for equal point masses and equal-mass $K_{16}$ strut networks, outer--inner vertex opposition ($V_0\overline{\mathbf p}_{\mathrm{outer}}=-V_1\overline{\mathbf p}_{\mathrm{inner}}$), non-overlapping papercraft nets with $N_{\mathrm{tabs}}=f_0-1$ glue tabs, and mirror-twin facet inversion ($\mathcal G_\lambda(r^2\overline S)=-\mathcal G_\lambda(S)$).

When a nonlinear geometric operation breaks this universal balance, it typically creates a three-component residual imbalance vector $\mathbf c=(c_x,c_y,c_z)\in\R^3$ that a single scalar parameter cannot zero out. In Dudeney's Group~VI ($304$ squares satisfying $\overline S=hS$ or $vS$), however, the $180^\circ$ proper half-turn symmetry ($C_2\subset SO(3)$) forces two coordinates of every residual vector to vanish identically, locking opposing components onto a single coordinate axis. This axial alignment explains why one-parameter restoration occurs throughout Group~VI: the eight planar Pell-minimum squares (\#70, $v/u=9/4\approx\sqrt{5}$), the integer power-law spoke roots ($\alpha^*=-1$ on \#45 and $\alpha^*=+3$ on \#286), all $107$ 3D power-law spoke roots $\beta^*\ne1$, the three non-associative balanced solids (\#86, \#739, \#756), all $35$ two-density ``diamond-in-a-diamond'' solids ($33$ denser cores led by \#673 plus the two lighter/hollow-cavity cores \#123 and \#54), all $158$ balanced nested wireframes ($d_1/d_0$), all $80$ balanced nested shells ($\sigma_1/\sigma_0$), and the chiral knight's-move mirror-twin towers (\#763).

Finally, when a nonlinear observable requires unweighted equilibrium without an opposing parameter—such as equal-weight unit-sphere spokes ($\mathbf V_0=\mathbf 0$), full 3D cubic spoke vanishing ($\mathbf V_3=\mathbf 0$), 3D entry-weighted rotational axisymmetry ($Q_{xy}=0$ and $\mathbf V_3=\mathbf 0$), or uniform-diameter wireframes and hollow shells—balance holds if and only if central inversion $\mathbf p_{-c}(\lambda)=-\mathbf p_c(\lambda)$ cancels opposite cells pairwise, singling out the $48$ associative squares of Dudeney's Group~III ($r^2S=\overline S$) led by D\"urer's \#175.

\begin{table}[htbp]
\centering
\caption{Synthesis of the four progressive mechanical tests across the $880$ order-four magic squares, organized by Baseline Equilibrium, Nonlinear Symmetry Breaking, and Dual Restoration (certified in exact rational arithmetic; see \Cref{sec:appendix}).}
\label{tab:synthesis}
\scriptsize
\setlength{\tabcolsep}{3.0pt}
\resizebox{\textwidth}{!}{%
\begin{tabular}{@{}l c l l@{}}
\toprule
Mechanical Regime & $n=4$ Count & Universal Law, Breakdown, or Restoration Threshold & Flagship Exemplars\\
\midrule
\multicolumn{4}{@{}l}{\textbf{1. Radial Spokes from the Origin (\Cref{sec:vectors})}}\\
Linear 3D \& 2D spokes ($\beta=1,\alpha=1$) & $880/880$ ($100\%$) & Baseline: $\sum\mathbf p_{ij}(\lambda)=\mathbf 0$, $\mathbf B_{z,r}=\mathbf 0$ identically & All semi-magic\\
Planar unit spokes ($\alpha=0$, $r_2/r_1=\sqrt{5}$) & $0/880$ vanish & Breakdown: $\min\|\mathbf B_{\mathrm{unit}}\|_2=\frac{18\sqrt{10}}{5}-8\sqrt{2}\approx0.0705$ ($v/u=9/4$) & $8$ in VI (\#70, \#76, \dots)\\
Planar power-law restoration ($\ell=r^{\alpha^*}$) & $204/880$ ($23.2\%$) & Parallel ring dipoles $\mathbf d_{\mathrm{out}}=\kappa\mathbf d_{\mathrm{in}}$ ($25$ at $\alpha^*=-1$, $25$ at $+3$) & \#45 ($-1$), \#286 ($+3$)\\
3D unit sphere $\sum\hat{\mathbf p}_{ij}(\lambda)=\mathbf 0$ ($\beta=0$) & $48/880$ ($5.5\%$) & Pointwise antipodal cancellation ($\iff$ Group~III); $107$ VI have $\beta^*\ne1$ & \#175 ($\mathbf 0$), \#70 ($\beta^*\approx2.24$)\\
3D unit-sphere tip weights $\mathbf w^*>\mathbf 0$ & $880/880$ ($100\%$) & Restoration: $\sum w_{ij}^*\hat{\mathbf p}_{ij}=\mathbf 0$, $\|\mathbf w^*-\mathbf 1\|_2\in[0.005,1.60]$ & All $880$ squares\\
3D cubic spoke law ($\beta=3$, $L_{ij}=R_{ij}^3$) & $880/880$ ($100\%$) & Universal vertical vanishing $(\mathbf V_3)_z=0$; $\mathbf V_3=\mathbf 0\iff$ Group~III ($48$) & All semi-magic ($z=0$)\\
\midrule
\multicolumn{4}{@{}l}{\textbf{2. Rotational Dynamics (\Cref{sec:dynamics})}}\\
Equal-mass point rotor $I_\perp(\lambda)/I_{zz}$ & $880/880$ ($100\%$) & Baseline: $\frac{1+17\lambda^2}{2}$ (spherical top $1$ at $\lambda=\frac{1}{\sqrt{17}}$, prolate $9$ at $\lambda=1$) & All semi-magic\\
Entry-weighted ($m_{ij}=s_{ij}$) Dzhanibekov split & $824/880$ flip & Breakdown: $I_{1,2}=170\mp|Q_{xy}|$ ($56$ have $Q_{xy}=0$; $48$ Group~III in 3D) & \#46 ($|Q_{xy}|=60$), \#175 ($0$)\\
\midrule
\multicolumn{4}{@{}l}{\textbf{3. Convex Hulls, Skeletons, and Shells (\Cref{sec:nested})}}\\
Uniform solid hull $\mathcal G_0(S)$ ($\mathbf c_0=\mathbf 0$) & $51/880$ ($5.8\%$) & Breakdown on $829$ gems; $48$ Group~III $+ 3$ Group~VI balance naturally & \#175, \#86, \#739, \#756\\
Outer--inner vertex opposition & $862/880$ ($V_0<16$) & Restoration law: $V_0\overline{\mathbf p}_{\mathrm{outer}}+V_1\overline{\mathbf p}_{\mathrm{inner}}=\mathbf 0$ ($511$ have 3D core $\mathcal G_1$) & All $V_0<16$ gems\\
Two-density solid $(\mathcal G_0,\mathcal G_1)$ & $35/511$ 3D cores & Group~VI: $33$ denser ($8$ with $\Delta\rho_1/\rho_0\le\frac{40}{19}$) $+ 2$ lighter/hollow cavity & \#673 ($\frac54$), \#45$\leftrightarrow$\#54\\
Non-negative inner-core point weights & $528/880$ ($60.0\%$) & $51$ natural $+ 477$ with $-\mathbf c_0\in\cone(P_S^{(1)})$; bead inlay balances $880/880$ & $528$ magic gems\\
Complete $K_{16}$ strut network ($120$ rods) & $880/880$ ($=\text{mass}$) & $\mathbf c=\mathbf 0$, $\Sigma_{\mathrm{rods}}=\frac{29}{45}\Sigma(\lambda)$; uniform-diam.\ drift cut $3.6\times$ at $\lambda=0.20$ & All $880$ (\#175)\\
Two-gauge nested wireframe $G(\mathcal G_0)\cup G(\mathcal G_1)$ & $206/880$ ($\lambda=1$) & $48$ Group~III $+ 158$ Group~VI ($164$ at $\lambda=0.20$) via rod ratio $d_1/d_0$ & \#673 ($d_1/d_0=0.793$)\\
Papercraft nets \& nested shells $\partial\mathcal G_0$ & $880/880$ nets & $100\%$ non-overlapping, $N_{\mathrm{tabs}}=f_0-1\in[7,15]$; $48+80$ VI shells balance & \#175 ($11$ tabs, $13$ folds)\\
\midrule
\multicolumn{4}{@{}l}{\textbf{4. Freestanding Gravity and Mirror-Twin Towers (\Cref{sec:stacking})}}\\
Single-block facet resting stability & $11{,}044/14{,}824$ & Cubic aspect $\lambda=0.20$ raises stability from $43.0\%$ ($\lambda=1$) to $74.5\%$ & $\max\theta_{\mathrm{tip}}=47.3^\circ$\\
Mirror-twin vertical towers $\mathcal G_\lambda(r^2\overline S)$ & $48/880$ ($5.5\%$) & Central inversion $-\mathcal G_\lambda(S)$ locks slanted facets ($17$ area-$5$ knight's-move) & \#763$\leftrightarrow$\#756, \#603\\
\bottomrule
\end{tabular}%
}
\end{table}

These findings suggest several natural directions for further investigation. In the planar spoke regime, the twenty-four non-central cells of $A_5$ at order $n=5$ lie on five concentric rings of radii $r\in\{1,\sqrt{2},2,\sqrt{5},2\sqrt{2}\}$, coupling the independent square-free radicals $\sqrt{2}$ and $\sqrt{5}$; it remains open whether the unit-spoke resultant $\mathbf B_{\mathrm{unit}}(S)=\sum z_{ij}\mathbf u_{ij}$ can ever vanish at order $n\ge5$, or how rapidly $\min_{S\in\mathcal M_n}\|\mathbf B_{\mathrm{unit}}(S)\|_2$ decays as $n\to\infty$ via simultaneous Diophantine approximation. In the rotational and solid regimes, while point-weighted Magic Gems become exact spherical tops ($I_{xx}=I_{yy}=I_{zz}$) at $\lambda_{\mathrm{iso}}=(n^2+1)^{-1/2}$, the uniform polyhedral solid $\mathcal G_\lambda(S)$ retains a small transverse splitting $(I_3-I_2)/I_2\in[0.00004,0.0243]$ at $n=4$, raising the question of whether a magic square of some order $n\ge5$ (or a two-density composite $(\mathcal G_0,\mathcal G_1)$ at $n=4$) has a continuous solid inertia tensor that is an exact scalar multiple of the identity $I_3$. In the convex-peeling regime, the outer hull $\mathcal G_0(S_n)$ of the Siamese magic square of odd order $n\ge5$ strips off only twelve vertices \cite[Thm.~6.4]{mathewson2026gems}, leaving $n^2-12$ points on five sheared planar polygons whose asymptotic onion-peeling depth $d_{\mathrm{onion}}(S_n)$ and inner layer $f$-vectors remain to be determined. Finally, beyond vertical mirror-twin columns (\Cref{thm:mirror-tower}), it is natural to ask whether a Magic Gem $\mathcal G_\lambda(S)$ and its $D_4\times\langle\tau\rangle$ images can tile a horizontal slab or interlock in three-dimensional space with full facet-to-facet contact.

\section*{Acknowledgments}
\label{sec:ai_statement}

This work grew out of the author's physical and computational explorations of 3D-printed Magic Gems, ball-and-rod wireframes, and folded papercraft models \cite{mathewson2025magicgemweb}, and forms the third part of a series alongside \cite{mathewson2026gems,mathewson2026moments}. Large language models (OpenAI's ChatGPT, Anthropic's Claude, and Google's Gemini) were used by the author as interactive assistants for symbolic checks, verification scripting, interactive simulator prototyping, and manuscript copy-editing; all mathematical formulations, proofs, and exact computational certificates were verified by the author, who takes full responsibility for the work.

\bibliographystyle{plain}
{\small\setlength{\itemsep}{0pt}\bibliography{references}}

\appendix
\section{Exact computational methods and certificates}
\label{sec:appendix}

All computations in this paper use Python~3.9 with NumPy~2.0, SciPy~1.13, and Python's exact rational module \texttt{fractions.Fraction}. The complete verification script \texttt{anc/verify\_paper\_c.py} and figure generator \texttt{anc/make\_figures\_c.py} run in under $10$\,s on a single core and write the certified census archive \texttt{anc/results\_paper\_c.json}.

\paragraph{Exact simplex moments for continuous polyhedra.}
As in \cite[App.~A]{mathewson2026gems}, doubling the centred coordinates \eqref{eq:coords} produces integer vectors $P=2\mathbf p\in\Z^3$, and every facet plane $aX+bY+cZ+d=0$ of $\mathcal G_0(S)$ and of the inner core gem $\mathcal G_1(S)=\conv P_S^{(1)}$ is certified in exact integer arithmetic. To compute the exact rational volume $\vol(\mathcal G)$, solid centroid $\mathbf c_{\mathrm{solid}}(\mathcal G)$, and second-moment covariance matrix $\Sigma_{\mathrm{solid}}(\mathcal G)$ of any 3-dimensional gem $\mathcal G$, we decompose $\mathcal G$ into tetrahedra $\Delta_k=\conv\{\mathbf v_0,\mathbf v_{k,1},\mathbf v_{k,2},\mathbf v_{k,3}\}$ sharing the rational vertex $\mathbf v_0=\frac12 P_0\in(\frac12\Z)^3$. Each tetrahedron has exact rational volume
\[
V_k = \frac{1}{6}\,\bigl|\det(\mathbf v_{k,1}-\mathbf v_0,\, \mathbf v_{k,2}-\mathbf v_0,\, \mathbf v_{k,3}-\mathbf v_0)\bigr| \in \Q_{\ge0},
\]
exact first moment
\[
\iiint_{\Delta_k} \mathbf x\,dV = \frac{V_k}{4}\,\bigl(\mathbf v_0 + \mathbf v_{k,1} + \mathbf v_{k,2} + \mathbf v_{k,3}\bigr) \in \Q^3,
\]
and exact canonical second moment
\[
\iiint_{\Delta_k} x_i x_j\,dV = \frac{V_k}{20}\left(2\sum_{a=0}^3 v_{a,i}v_{a,j} + \sum_{0\le a\ne b\le3} v_{a,i}v_{b,j}\right) \in \Q,
\]
where $(\mathbf v_{k,0},\mathbf v_{k,1},\mathbf v_{k,2},\mathbf v_{k,3})=(\mathbf v_0,\mathbf v_{k,1},\mathbf v_{k,2},\mathbf v_{k,3})$. Summing over all simplices $\Delta_k$ yields $\vol(\mathcal G)=\sum_k V_k\in\Q$, $\mathbf c_{\mathrm{solid}}(\mathcal G)=\frac{1}{\vol(\mathcal G)}\sum_k\iiint_{\Delta_k}\mathbf x\,dV\in\Q^3$, and $\Sigma_{\mathrm{solid},ij}=\frac{1}{\vol(\mathcal G)}\sum_k\iiint_{\Delta_k}x_ix_j\,dV - c_{\mathrm{solid},i}c_{\mathrm{solid},j}\in\Q$ without floating-point rounding. Because the simplicial decomposition partitions the convex body $\mathcal G$, these rational integrals are independent of the choice of boundary triangulation.

\paragraph{Complete list of the 33 two-density balanced gems.}
For completeness, \Cref{tab:all33} lists all $33$ order-four magic squares for which $\dim\mathcal G_1(S)=3$ and the solid centroids $\mathbf c_0(S)$ and $\mathbf c_1(S)$ are strictly anti-parallel, together with their exact rational balancing density increment $\Delta\rho_1/\rho_0=-\vol(\mathcal G_0)c_0/(\vol(\mathcal G_1)c_1)>0$ (\Cref{thm:diamond-balance}). Every one of the $33$ squares belongs to Dudeney's Group~VI.

\begin{table}[htbp]
\centering
\caption{Complete census of the $33$ order-four Magic Gems admitting exact two-density ``diamond-in-a-diamond'' solid balance $\mathbf c_{\mathrm{comp}}=\mathbf 0$ (\Cref{thm:diamond-balance}), sorted by increasing excess density ratio $\Delta\rho_1/\rho_0$.}
\label{tab:all33}
\scriptsize
\setlength{\tabcolsep}{4.5pt}
\begin{tabular}{r c c c | r c c c | r c c c}
\toprule
Gem & $\vol(\mathcal G_0)$ & $\vol(\mathcal G_1)$ & $\Delta\rho_1/\rho_0$ &
Gem & $\vol(\mathcal G_0)$ & $\vol(\mathcal G_1)$ & $\Delta\rho_1/\rho_0$ &
Gem & $\vol(\mathcal G_0)$ & $\vol(\mathcal G_1)$ & $\Delta\rho_1/\rho_0$\\
\midrule
\#673 & $257/3$ & $56/3$ & $5/4$   & \#665 & $281/3$ & $6$    & $61/13$  & \#697 & $248/3$ & $28/3$ & $166/11$\\
\#502 & $271/3$ & $6$    & $4/3$   & \#252 & $233/3$ & $7$    & $88/9$   & \#671 & $84$    & $2/3$  & $31/2$\\
\#45  & $214/3$ & $6$    & $85/61$ & \#166 & $93$    & $8/3$  & $21/2$   & \#383 & $284/3$ & $7$    & $47/3$\\
\#50  & $254/3$ & $22/3$ & $28/19$ & \#158 & $88$    & $9$    & $57/5$   & \#398 & $260/3$ & $7$    & $95/6$\\
\#307 & $277/3$ & $22/3$ & $23/15$ & \#168 & $280/3$ & $7$    & $110/9$  & \#18  & $71$    & $9$    & $33/2$\\
\#49  & $81$    & $2$    & $2$     & \#672 & $260/3$ & $4/3$  & $99/8$   & \#627 & $83$    & $4/3$  & $67/4$\\
\#38  & $263/3$ & $22/3$ & $31/15$ & \#167 & $90$    & $8/3$  & $27/2$   & \#452 & $236/3$ & $6$    & $118/7$\\
\#47  & $247/3$ & $22/3$ & $40/19$ & \#199 & $256/3$ & $22/3$ & $206/15$ & \#604 & $221/3$ & $26/3$ & $229/10$\\
\#261 & $259/3$ & $7$    & $47/12$ & \#615 & $244/3$ & $5$    & $14$     & \#226 & $235/3$ & $7$    & $143/6$\\
\#670 & $244/3$ & $10$   & $103/25$& \#6   & $239/3$ & $7$    & $127/9$  & \#626 & $69$    & $2/3$  & $49/2$\\
\#488 & $277/3$ & $22/3$ & $71/16$ & \#573 & $251/3$ & $8$    & $44/3$   & \#436 & $235/3$ & $6$    & $191/6$\\
\bottomrule
\end{tabular}
\end{table}

\end{document}